\documentclass{amsart}
\usepackage{amssymb}
\usepackage{epsfig}
\newtheorem{theorem}{Theorem}[section]

\newtheorem{proposition}[theorem]{Proposition}

\newtheorem{definition}[theorem]{Definition}
\newtheorem{conj}[theorem]{Conjecture}

\newcommand{\C}{\mathbb{C}}
\newcommand{\Z}{\mathbb{Z}}
\newcommand{\CP}{\mathbb{CP}}
\newcommand{\R}{\mathbb{R}}
\newcommand{\Q}{\mathbb{Q}}
\newcommand{\dbar}{\bar{\partial}}
\newcommand{\comment}[1]{}
\title{Monodromy invariants in symplectic topology}
\author{Denis Auroux}
\address{Centre de Math\'ematiques, Ecole Polytechnique, F-91128 Palaiseau,
France}
\email{auroux@math.polytechnique.fr}

\begin{document}

\maketitle

%\begin{center}
%\bf \fbox{NOTES FOR LECTURE 3}
%\end{center}
%\bigskip

\textit{
This text is a set of lecture notes for a series of four talks
given at I.P.A.M., Los Angeles, on March 18-20, 2003. 
The first lecture provides a quick overview of symplectic topology and
its main tools: symplectic manifolds, almost-complex structures,
pseudo-holomorphic curves, Gromov-Witten invariants and Floer homology. The
second and third lectures focus on symplectic Lefschetz pencils: existence
(following Donaldson), monodromy, and applications to symplectic topology,
in particular the connection to Gromov-Witten invariants of symplectic
4-manifolds (following Smith) and to Fukaya categories (following Seidel).
In the last lecture, we offer
an alternative description of symplectic 4-manifolds by viewing them as
branched covers of the complex projective plane; the corresponding monodromy
invariants and their potential applications are discussed.
}

\section{Introduction to symplectic topology}

%\comment{

In this lecture, we recall basic notions about symplectic manifolds, and
briefly review some of the main tools used to study them. Most of the
topics discussed here can be found in greater detail in standard graduate
books such as \cite{McS}.

\subsection{Symplectic manifolds}

\begin{definition}
A {\em symplectic structure} on a smooth manifold $M$ is a closed
non-degenerate $2$-form $\omega$, i.e.\ an element $\omega\in\Omega^2(M)$
such that $d\omega=0$ and $\forall v\in TM-\{0\}$, $\iota_v\omega\neq 0$.
\end{definition}

For example, $\R^{2n}$ carries a standard symplectic structure, given by
the 2-form $\omega_0=\sum dx_i\wedge dy_i$. Similarly, every orientable
surface is symplectic, taking for $\omega$ any non-vanishing volume form.

Since $\omega$ induces a non-degenerate antisymmetric bilinear pairing on
the tangent spaces to $M$, it is clear that every symplectic manifold is
even-dimensional and orientable (if $\dim M=2n$, then $\frac{1}{n!}\omega^n$
defines a volume form on $M$).

Two important features of symplectic structures that set them apart from
most other geometric structures are the existence of a large number of
symplectic automorphisms, and the absence of local geometric invariants.

The first point is illustrated by the following construction. Consider
a smooth function $H:M\to\R$ (a {\it Hamiltonian}), and define $X_H$ to be
the vector field on $M$ such that $\omega(X_H,\cdot)=dH$. Let
$\phi_t:M\to M$ be the family of diifeomorphisms generated by the flow of
$X_H$, i.e., $\phi_0=\mathrm{Id}$ and $\frac{d}{dt}\phi_t(x)=X_H(\phi_t(x))$.
Then $\phi_t$ is a {\it symplectomorphism}, i.e.\ $\phi_t^*\omega=\omega$.
Indeed, we have $\phi_0^*\omega=\omega$, and
$$\frac{d}{dt}\phi_t^*\omega=\phi_t^*(L_{X_H}\omega)=
\phi_t^*(d\iota_{X_H}\omega+\iota_{X_H}d\omega)=\phi_t^*(d(dH)+0)=0.$$
Therefore, the group of symplectomorphisms $\mathrm{Symp}(M,\omega)$ is
infinite-dimensional, and its Lie algebra contains all Hamiltonian vector
fields. This is in contrast with the case of Riemannian metrics, where
isometry groups are much smaller.

The lack of local geometric invariants of symplectic structures is 
illustrated by two classical results of fundamental importance, which show
that the study of symplectic manifolds largely reduces to {\it topology}
(i.e., to discrete invariants): Darboux's theorem, and Moser's stability 
theorem. The first one shows that all symplectic forms are locally equivalent,
in sharp contrast to the case of Riemannian metrics where curvature provides
a local invariant, and the second one shows that exact deformations of
symplectic structures are trivial.

\begin{theorem}[Darboux]
Every point in a symplectic manifold $(M^{2n},\omega)$ admits a neighborhood
that is symplectomorphic to a neighborhood of the origin in
$(\R^{2n},\omega_0)$.
\end{theorem}

\proof
We first use local coordinates to map a neighborhood of a given
point in $M$ diffeomorphically onto a neighborhood $V$ of the origin in
$\R^{2n}$. Composing this diffeomorphism $f$ with a suitable linear
transformation of $\R^{2n}$, we can ensure that the symplectic form
$\omega_1=(f^{-1})^*\omega$ coincides with $\omega_0$ at the origin.
This implies that, restricting to a smaller neighborhood if necessary,
the closed 2-forms $\omega_t=t\,\omega_1 + (1-t)\omega_0$ are non-degenerate
over $V$ for all $t\in [0,1]$.

Using the Poincar\'e lemma, consider a family of $1$-forms $\alpha_t$ on $V$
such that $\frac{d}{dt}\omega_t=-d\alpha_t$. Subtracting a constant $1$-form
from $\alpha_t$ if necessary, we can assume that $\alpha_t$ vanishes at the
origin for all $t$. Using the non-degeneracy of $\omega_t$
we can find vector fields $X_t$ such that $\iota_{X_t}\omega_t = \alpha_t$.
Let $(\phi_t)_{t\in [0,1]}$ be the flow generated by $X_t$, i.e.\ the family
of diffeomorphisms defined by $\phi_0=\mathrm{Id}$, $\frac{d}{dt}\phi_t(x)=
X_t(\phi_t(x))$; we may need to restrict to a smaller neighborhood
$V'\subset V$ of the origin in order to make the flow $\phi_t$ well-defined
for all $t$. We then have
$$\frac{d}{dt}\phi_t^*\omega_t=\phi_t^*(L_{X_t}\omega_t) +
\phi_t^*\bigl(\frac{d\omega_t}{dt}\bigr)=
\phi_t^*(d(\iota_{X_t}\omega_t)-d\alpha_t)=0,$$
and therefore $\phi_1^*\omega_1=\omega_0$. Therefore, $\phi_1^{-1}\circ f$
induces a symplectomorphism from a neighborhood of $x$ in $(M,\omega)$ to a
neighborhood of the origin in $(\R^{2n},\omega_0)$.
\endproof

\begin{theorem}[Moser]
Let $(\omega_t)_{t\in [0,1]}$ be a continuous family of symplectic forms
on a compact manifold $M$. Assume that the cohomology class $[\omega_t]\in
H^2(M,\R)$ does not depend on $t$. Then $(M,\omega_0)$ is symplectomorphic
to $(M,\omega_1)$.
\end{theorem}

\proof
We use the same argument as above: since $[\omega_t]$ is constant there exist
\hbox{1-forms}
$\alpha_t$ such that $\frac{d}{dt}\omega_t=-d\alpha_t$. Define vector fields
$X_t$ such that $\iota_{X_t}\omega_t=\alpha_t$ and the corresponding flow
$\phi_t$. By the same calculation as above, we
conclude that $\phi_1^*\omega_1 = \omega_0$.
\endproof

\subsection{Submanifolds in symplectic manifolds}

\begin{definition} A submanifold $W\subset (M^{2n},\omega)$ is called
{\em symplectic} if $\omega_{|W}$ is non-degenerate at every point of $W$
(it is then a symplectic form on $W$);
{\em isotropic} if $\omega_{|W}=0$; and {\em Lagrangian} if it is isotropic
of maximal dimension $\dim W=n=\frac{1}{2}\dim M$.
\end{definition}

An important example is the following: given any smooth manifold $N$, the
cotangent bundle $T^*N$ admits a canonical symplectic structure that can be
expressed locally as $\omega=\sum dp_i\wedge dq_i$ (where $(q_i)$ are local
coordinates on $N$ and $(p_i)$ are the dual coordinates on the cotangent
spaces). Then the zero section is a Lagrangian submanifold of $T^*N$.

Since the symplectic form induces a non-degenerate pairing between tangent
and normal spaces to a Lagrangian submanifold, the normal bundle to a
Lagrangian submanifold is always isomorphic to its cotangent bundle. The
fact that this isomorphism extends beyond the infinitesimal level is a
classical result of Weinstein:

\begin{theorem}[Weinstein]
For any Lagrangian submanifold $L\subset (M^{2n},\omega)$, there exists
a neighborhood of $L$ which is symplectomorphic to a neighborhood of the
zero section in the cotangent bundle $(T^*L,\sum dp_i\wedge dq_i)$.
\end{theorem}

There is also a neighborhood theorem for symplectic submanifolds; in that
case, the local model for a neighborhood of the submanifold $W\subset M$
is a neighborhood of the zero section in the symplectic vector bundle
$NW$ over $W$ (since $Sp(2n)$ retracts onto $U(n)$, the classification of
symplectic vector bundles is the same as that of complex vector bundles).

\subsection{Almost-complex structures}

\begin{definition}
An {\em almost-complex structure} on a manifold $M$ is an endomorphism $J$ of
the tangent bundle $TM$ such that $J^2=-\mathrm{Id}$. An almost-complex
structure $J$ is said to be {\em tamed} by a symplectic form $\omega$ if
for every non-zero tangent vector $u$ we have $\omega(u,Ju)>0$; it is
{\em compatible} with $\omega$ if it is $\omega$-tame and 
$\omega(u,Jv)=-\omega(Ju,v)$; equivalently, $J$ is $\omega$-compatible if
and only if $g(u,v)=\omega(u,Jv)$ is a Riemannian metric.
\end{definition}

\begin{proposition}
Every symplectic manifold $(M,\omega)$ admits a compatible almost-complex
structure. Moreover, the space of $\omega$-compatible (resp.\ $\omega$-tame)
almost-complex structures is contractible.
\end{proposition}

This result follows from the fact that the space of compatible (or tame)
complex structures on a symplectic vector space is non-empty and
contractible (this can be seen by constructing explicit retractions);
it is then enough to observe that a compatible
(resp.\ tame) almost-complex structure on a symplectic manifold is simply
a section of the bundle $End(TM)$ that defines a compatible (resp.\ tame)
complex structure on each tangent space.

An almost-complex structure induces a splitting of the complexified tangent
and cotangent bundles: $TM\otimes\C=TM^{1,0}\oplus TM^{0,1}$, where
$TM^{1,0}$ and $TM^{0,1}$ are respectively the $+i$ and $-i$ eigenspaces
of $J$ (i.e., $TM^{1,0}=\{v-iJv,\ v\in TM\}$, and similarly for $TM^{0,1}$;
for example, on $\C^n$ equipped with its standard complex structure, the
$(1,0)$ tangent space is generated by $\partial/\partial z_i$ and the
$(0,1)$ tangent space by $\partial/\partial \bar{z}_i$.
Similarly, $J$ induces a complex structure on the cotangent bundle, and
$T^*M\otimes\C=T^*M^{1,0}\oplus T^*M^{0,1}$ (by definition $(1,0)$-forms
are those which pair trivially with $(0,1)$-vectors, and vice versa).
This splitting of the cotangent bundle induces a splitting of differential
forms into ``types'': $\bigwedge^rT^*M\otimes \C=
\bigoplus_{p+q=r}\bigwedge^pT^*M^{1,0} \otimes \bigwedge^qT^*M^{0,1}$.
Moreover, given a function $f:M\to\C$ we can write $df=\partial f + \dbar
f$, where $\partial f=\frac{1}{2}(df-i\,df\circ J)$ and $\dbar f=\frac{1}{2}
(df+i\,df\circ J)$ are the $(1,0)$ and $(0,1)$ parts of $df$ respectively.
Similarly, given a complex vector bundle $E$ over $M$ equipped with a
connection, the covariant derivative $\nabla$ can be split into operators
$\partial^\nabla:\Gamma(E)\to\Omega^{1,0}(E)$ and
$\dbar{}^\nabla:\Gamma(E)\to\Omega^{0,1}(E)$.

Although the tangent space to a symplectic manifold $(M,\omega)$ equipped
with a compatible almost-complex structure $J$ can be pointwise identified
with $(\C^n,\omega_0,i)$, there is an important difference between a
symplectic manifold equipped with a compatible almost-complex structure and
a complex K\"ahler manifold: the possible lack of {\it integrability} of
the almost-complex structure, namely the fact that the Lie bracket of two
$(1,0)$ vector fields is not necessarily of type $(1,0)$.

\begin{definition}
The {\em Nijenhuis tensor} of an almost-complex manifold $(M,J)$ is the
quantity defined by $N_J(X,Y)=\frac{1}{4}([X,Y]+J[X,JY]+J[JX,Y]-[JX,JY])$.
The almost-complex structure $J$ is said to be {\em integrable} if $N_J=0$.
\end{definition}

It can be checked that $N_J$ is a tensor (i.e., only depends on the values
of the vector fields $X$ and $Y$ at a given point), and that
$N_J(X,Y)=2\,\mathrm{Re}([X^{1,0},Y^{1,0}]^{(0,1)})$. The non-vanishing
of $N_J$ is therefore an obstruction to the integrability of a local frame
of $(1,0)$ tangent vectors, i.e.\ to the existence of local holomorphic
coordinates. The Nijenhuis tensor is also related to the fact that the
exterior differential of a $(1,0)$-form may have a non-zero component of
type $(0,2)$, so that the $\partial$ and $\dbar$ operators on differential
forms do not have square zero ($\dbar{}^2$ can be expressed in terms of
$\partial$ and the Nijenhuis tensor).

\begin{theorem}[Newlander-Nirenberg] Given an almost-complex manifold
$(M,J)$, the following properties are
equivalent: $(i)$~$N_J=0$; $(ii)$ $[T^{1,0}M,T^{1,0}M]\subset T^{1,0}M$;
$(iii)$ $\dbar{}^2=0$; $(iv)$ $(M,J)$ is a complex manifold, i.e.\ admits
complex analytic coordinate charts.
\end{theorem}

\subsection{Pseudo-holomorphic curves and Gromov-Witten invariants}
\leavevmode\medskip

\noindent
Pseudo-holomorphic curves, first introduced by Gromov in 1985 \cite{Gr},
have since then become the most important tool in modern symplectic
topology. In the same way as the study of complex curves in complex
manifolds plays a central role in algebraic geometry, the study of
pseudo-holomorphic curves has revolutionized our
understanding of symplectic manifolds.

The equation for holomorphic maps between two almost-complex manifolds
becomes overdetermined as soon as the complex dimension of the domain
exceeds $1$, so we cannot expect the presence of any almost-complex
submanifolds of complex dimension $\ge 2$ in a symplectic manifold equipped with
a generic almost-complex structure. On the other hand, $J$-holomorphic
curves, i.e.\ maps from a Riemann surface $(\Sigma,j)$ to the manifold
$(M,J)$ such that $J\circ df = df\circ j$, are governed by an elliptic
PDE, and their study makes sense even in non-K\"ahler symplectic manifolds.
The questions that we would like to answer are of the following type:
\medskip

{\it Given a compact symplectic manifold $(M,\omega)$ equipped with a generic
compatible almost-complex structure $J$ and a homology class $\beta\in H_2(M,\Z)$,
what is the number of pseudo-holomorphic curves of given genus $g$,
representing the homology class $\beta$ and passing through $r$ given points
in $M$ (or through $r$ given submanifolds in $M$)~?}%
\medskip

The answer to this question is given by {\it Gromov-Witten invariants},
which count such curves (in a sense that is not always obvious, as the
result can e.g.\ be negative, and need not even be an integer). One starts
by introducing a {\it moduli space} $\mathcal{M}_{g,r}(\beta)$ of genus $g$
pseudo-holomorphic curves with $r$ marked points representing a given
homology class $\beta\in H_2(M,\Z)$. This moduli space comes equipped
with $r$ {\it evaluation maps} $ev_i:\mathcal{M}_{g,r}(\beta)\to M,$
associating to a curve the image of its $i$-th marked point. Considering
$r$ submanifolds of $M$ representing homology classes Poincar\'e dual to
cohomology classes $\alpha_1,\dots,\alpha_r\in H^*(M,\R)$, and assuming the
dimensions to be the right ones, we then want to define a
number \begin{equation}\mathrm{GW}_{g,\beta}(M;\alpha_1,\dots,\alpha_r)=
\int_{\mathcal{M}_{g,r}(\beta)} ev_1^*\alpha_1\wedge\dots\wedge ev_r^*
\alpha_r.\end{equation}
We can also choose to impose restrictions on the complex structure of
the domain (or on the positions of the marked points in the domain) by
introducing into the integral a cohomology class pulled back from the moduli
space of genus $g$ Riemann surfaces with $r$ marked points.

Giving a precise meaning to the right-hand side of (1) is a very delicate
task, which has been a central preoccupation of symplectic geometers
for more than ten years following Gromov's seminal work. We only give a
very simplified description of the main issues; the reader is
referred to the book \cite{McS2} for a detailed
discussion of the so-called ``weakly monotone'' case, and to more
recent work (Fukaya-Ono, Li-Tian \cite{LT}, Ruan, Siebert, Hofer-Salamon,
\dots) for the general case, which requires more sophisticated techniques. 

To start with, one must study deformations of
pseudo-holomorphic curves, by linearizing the equation $\dbar_J f=0$ near a
solution. The linearized Cauchy-Riemann operator $D_{\dbar}$, whose kernel
describes infinitesimal deformations of a given curve $(f:\Sigma\to M)\in
\mathcal{M}_{g,r}(\beta)$, is a Fredholm operator of (real) index
$$2d:=\mathrm{ind}\,D_{\dbar}=(1-g)(\dim_\R M-6)+2r+2\,c_1(TM)\cdot\beta.$$

When the considered curve is {\it regular}, i.e.\ when the linearized operator
$D_{\dbar}$ is surjective, the deformation theory is unobstructed,
and if the curve has no automorphisms (which is the case of a generic curve
provided that $r\ge 3$ when $g=0$ and $r\ge 1$ when $g=1$), we expect the
moduli space $\mathcal{M}_{g,r}(\beta)$ to be locally a smooth manifold of
real dimension $2d$.

The main result underlying the theory of pseudo-holomorphic curves is
Gromov's compactness theorem (see \cite{Gr}, \cite{McS2}, \dots):

\begin{theorem}[Gromov]
Let $f_n:(\Sigma_n,j_n)\to (M,\omega,J)$ be a sequence of pseudo-holomorphic
curves in a compact symplectic manifold, representing a fixed homology
class. Then a subsequence of $\{f_n\}$
converges (in the ``Gromov-Hausdorff topology'') to a limiting map
$f_\infty$, possibly singular.
\end{theorem}

The limiting curve $f_\infty$ can have a very complicated structure, and in
particular its domain may be a nodal Riemann surface with more than one
component, due to the phenomenon of {\it bubbling}. For example, the
sequence of degree $2$ holomorphic curves $f_n:\CP^1\to\CP^2$ defined by
$f_n(u\!:\!v)=(u^2\!:\!uv\!:\!\frac{1}{n}v^2)$ converges to a singular
curve with two degree 1 components: for $(u\!:\!v)\neq
(0\!:\!1)$, we have $\lim f_n(u\!:\!v)=(u\!:\!v\!:\!0)$, so that the
sequence apparently converges to a line in $\CP^2$. However the derivatives
of $f_n$ become unbounded near $(0\!:\!1)$, and composing $f_n$ with the
coordinate change $\phi_n(u\!:\!v)=(\frac{1}{n}u\!:\!v)$ we obtain
$f_n\circ \phi_n(u\!:\!v)=(\frac{1}{n^2}u^2\!:\!\frac{1}{n}uv\!:\!
\frac{1}{n}v^2)=(\frac{1}{n}u^2\!:\!uv\!:\!v^2)$, which
converges to $(0\!:\!u\!:\!v)$ everywhere except at
$(1\!:\!0)$, giving the other component (the ``bubble'') in the
limiting curve.

The presence of a symplectic structure on the target manifold $M$ is crucial
for compactness, as it provides an a priori estimate of the {\it energy} of
a pseudo-holomorphic curve, and hence makes it possible to control the
bubbling phenomenon: for a $J$-holomorphic map $f:(\Sigma,j)\to
(M,\omega,J)$, we have
$$\int_\Sigma \|df\|^2=\int_\Sigma f^*\omega=[\omega]\cdot f_*[\Sigma].$$

The definition of Gromov-Witten invariants requires a compactification of
the moduli space of pseudo-holomorphic curves, in order to be able to define
a fundamental cycle for $\mathcal{M}_{g,r}(\beta)$ against which cohomology
classes can be evaluated. It follows from Gromov's compactness theorem that
this compactification can be achieved by considering the moduli space
$\overline{\mathcal{M}}_{g,r}(\beta)$ of {\em stable maps}, i.e.\ 
$J$-holomorphic maps $f:\bigsqcup (\Sigma_\alpha,j_\alpha)\to M$ with domain
a tree of Riemann surfaces, and with a discrete set of automorphisms (which
imposes conditions on the genus $0$ or $1$ ``ghost components'').
If we assume that none of the stable maps in the compactified moduli space
has multiply covered components (i.e.\ components that factor through
non-trivial coverings of Riemann surfaces), then the moduli space of stable
maps can be used to define a fundamental class
$[\overline{\mathcal{M}}_{g,r}(\beta)]$. However, because the presence of
curves with non-trivial automorphisms makes the moduli space an orbifold
rather than a manifold, the fundamental class is in general a {\it rational}
homology class in $H_{2d}(\overline{\mathcal{M}}_{g,r}(\beta),\mathbb{Q})$.

Technically, the hardest case is when some of the stable maps in the moduli
space have components which are multiply covered. The stable map
compactification can then fail to provide a suitable fundamental cycle,
because the actual dimension of the moduli space of pseudo-holomorphic curves
may exceed that predicted by the index calculation, even for a generic choice of
$J$, due to the possibility of arbitrarily moving the branch points of the
multiple components. It is then necessary to break the symmetry and restore
transversality by perturbing the holomorphic curve equation $\dbar_J f=0$ 
into another equation $\dbar_J f=\nu(f)$ whose space of solutions
has the correct dimension (\cite{LT}, \dots), making it
possible to define a {\it virtual fundamental cycle}
$[\overline{\mathcal{M}}_{g,r}(\beta)]^{vir}$. For the same reasons as
above, this fundamental cycle is in general a {\it rational}
homology class in $H_{2d}(\overline{\mathcal{M}}_{g,r}(\beta),\mathbb{Q})$.

\subsection{Floer homology for Hamiltonians}
Floer homology has been introduced in symplectic geometry as an attempt to
prove the {\it Arnold conjecture} for fixed points of Hamiltonian
diffeomorphisms: for every non-degenerate Hamiltonian diffeomorphism
$\phi:(M,\omega)\to (M,\omega)$ (i.e., the flow of a time-dependent family
of Hamiltonian vector fields), the number of fixed points of $\phi$ is
bounded from below by the sum of the Betti numbers of $M$. Briefly speaking,
the idea is to obtain invariants from the Morse
theory of a functional defined on an infinite-dimensional space, and hence
to deduce existence results for critical points of the functional.

To define the Floer homology of the diffeomorphism of $(M,\omega)$ generated
by a 1-periodic family of Hamiltonians $H:S^1\times M\to\R$, one introduces
the {\it action functional} on the space of contractible loops in $M$: given
a loop $\gamma:S^1\to M$ bounding a disk $u:D^2\to M$, we define
$$\mathcal{A}_H(\gamma)=-\int_{D^2} u^*\omega - \int_{S^1}
H(t,\gamma(t))\,dt.$$
The critical points of $\mathcal{A}_H$ correspond to the 1-periodic closed
orbits of the Hamiltonian flow, and the gradient trajectories correspond to
maps $f:\R\times S^1\to M$ satisfying an equation of the form
\begin{equation}
\frac{\partial f}{\partial s}+J(f)\frac{\partial f}{\partial t}-\nabla
H(t,f)=0.\end{equation}
When $H$ is independent of $t$, and if we start with a constant
loop $\gamma$, a gradient trajectory of $\mathcal{A}_H$ is simply a gradient
trajectory of $H$; however, when $H=0$, the equation for gradient
trajectories reduces to that of pseudo-holomorphic curves.

The {\it Floer complex} $CF_*(M,H)$ is a free graded module with one generator
for each contractible 1-periodic orbit $\gamma:S^1\to M$ of the Hamiltonian
flow, or more precisely for each pair $(\gamma,[u])$, where $[u]\in
\pi_2(M,\gamma)$ is a relative homotopy class (taking this more refined
description, $CF_*(M,H)$ is naturally a module over the {\it Novikov ring}).
The grading is defined by
the {\em Conley-Zehnder index} of the critical point $(\gamma,[u])$ of
$\mathcal{A}_H$. Given two periodic orbits $\gamma^-$ and $\gamma^+$, we
can consider the moduli space $\mathcal{M}(\gamma^-,\gamma^+;A)$ of gradient
trajectories joining $\gamma^-$ to $\gamma^+$ (i.e., of solutions to the
gradient flow equation (2) that converge to $\gamma^-$ for $s\to -\infty$
and to $\gamma^+$ for $s\to +\infty$) realizing a given relative homology
class $A$. The expected dimension of this moduli space, which always carries
an $\R$-action by translation in the $s$ direction, is the difference betwen
the Conley-Zehnder indices of $(\gamma^-,[u])$ and $(\gamma^+,[u]\#A)$;
hence, assuming regularity and compactness of the moduli spaces
$\mathcal{M}(\gamma^-,\gamma^+;A)$ we can define an operator $\partial$ on
$CF_*(M,H)$ by the formula
$$\partial (\gamma^-,[u])=\sum_{\gamma^+,A}
\#(\mathcal{M}(\gamma^-,\gamma^+;A)/\R)\ (\gamma^+,[u]\#A),$$
where the sum runs over pairs $(\gamma^+,A)$ such that
the expected dimension of the moduli space is equal to $1$.
After a suitable discussion of regularity and compactness (bubbling) issues,
which are essentially identical to those encountered in the definition of
Gromov-Witten invariants, a proper meaning can be given to this definition
of the operator $\partial$, and the following result can be obtained:

\begin{theorem}
For regular $(H,J)$, the Floer differential satisfies
$\partial^2=0$, and hence we can define the {\it Floer homology}
$HF_*(M,\omega,H,J)=\mathrm{Ker}\,\partial/\mathrm{Im}\,\partial$.
Moreover, the Floer homologies obtained for different $(H,J)$ are
naturally isomorphic.
\end{theorem}

Moreover, by considering the limit as $H\to 0$ one can relate Floer
homology to either the quantum or classical homology of $M$ (with
coefficients in the Novikov ring), which yields an inequality between
the number of critical points of $\mathcal{A}_H$ and the total rank
of $H_*(M)$, and hence the Arnold conjecture.

The technically easier {\it monotone} case of these results is due to Floer;
the general
case has been treated subsequently using the same methods as for the
definition of general Gromov-Witten invariants (the reader is referred to
e.g.\ \cite{McS2} and \cite{Sa} for detailed expositions on Hamiltonian
Floer homology).

\subsection{Floer homology for Lagrangians}
Floer homology has also been successfully introduced for the study of
Lagrangian submanifolds, a construction which has taken a whole new
importance after Kontsevich's formulation of the homological mirror
symmetry conjecture.

Consider two compact orientable
(relatively spin) Lagrangian submanifolds $L_0$ and $L_1$ in
a symplectic manifold $(M,\omega)$ equipped with a compatible almost-complex
structure $J$. Lagrangian Floer homology corresponds to the Morse theory of
a functional on (a covering of) the space of arcs joining $L_0$ to $L_1$,
whose critical points are constant paths.

The {\it Floer complex} $CF^*(L_0,L_1)$ is
a free module with one generator for each intersection
point $p\in L_0\cap L_1$, and grading given by the
{\it Maslov index}. To be more precise, as in the case of Hamiltonian Floer
homology, in the general case one needs to consider pairs $(p,[u])$
where $[u]$ is the equivalence class of a map $u:[0,1]\times [0,1]\to M$ such
that $u(\cdot,0)\in L_0$, $u(\cdot,1)\in L_1$, $u(1,\cdot)=p$,
and $u(0,\cdot)$ is a fixed arc joining $L_0$ to $L_1$ (two maps $u,u'$ are
equivalent if they have the same symplectic area and the same Maslov index).

Given two points $p_{\pm}\in L_0\cap L_1$, we can define a moduli space
$\mathcal{M}(p_-,p_+;A)$ of pseudo-holomorphic maps $f:\R\times [0,1]\to M$ 
such that $f(\cdot,0)\in L_0$, $f(\cdot,1)\in L_1$, and $\lim_{t\to \pm\infty}
f(t,\tau)=p_{\pm}$ $\forall \tau\in [0,1]$, realizing a given relative homology class $A$;
the expected dimension of this moduli space is the difference of Maslov
indices. Assuming regularity and compactness of $\mathcal{M}(p_-,p_+;A)$,
we can define an operator $\partial$ on $CF^*(L_0,L_1)$ by the formula
$$\partial (p_-,[u])=\sum_{p_+,A}
\#(\mathcal{M}(p_-,p_+;A)/\R)\ (p_+,[u]\#A),$$
where the sum runs over pairs $(p_+,A)$ for which the expected dimension of
the moduli space is $1$.

In all good cases we expect to have $\partial^2=0$, which allows us to define
the Floer homology 
$HF^*(L_0,L_1)=\mathrm{Ker}\,\partial/\mathrm{Im}\,\partial$.
However, a serious technical difficulty arises from the fact that, in the case
of the moduli spaces $\mathcal{M}(p_-,p_+;A)$ of pseudo-holomorphic strips,
bubbling can occur on the boundary of the domain, so that limit curves can
contain both $S^2$ and $D^2$ bubble components. Because disc bubbling is a
phenomenon that arises in real codimension 1, the fundamental chain
associated to the compactified moduli space is not always a cycle; in that
case, to define Floer homology we need to add other curves to the moduli
space $\mathcal{M}(p_-,p_+;A)$ in order to obtain a cycle. This gives rise
to a series of obstructions that may prevent Floer homology from being
well-defined. The detailed analysis of the moduli space and of the
obstructions to the definition of Lagrangian Floer homology has been
carried out in the recent work of Fukaya, Oh, Ohta and Ono \cite{FO3},
refining the earlier work of Floer and Oh on technically easier special cases.

When Floer homology is well-defined, it has important consequences on the
intersection properties of Lagrangian submanifolds. Indeed, for every
{\it Hamiltonian} diffeomorphism $\phi$ we have 
$HF^*(L_0,L_1)=HF^*(L_0,\phi(L_1))$; and if $L_0$ and $L_1$ intersect
transversely, then the total rank of $HF^*(L_0,L_1)$ gives a lower bound on
the number of intersection points of $L_0$ and $L_1$. Moreover, $HF^*(L,L)$
is related to the usual cohomology $H^*(L)$ via a spectral sequence which
degenerates in the case where $i_*:H_*(L,\mathbb{Q})\to H_*(M,\mathbb{Q})$
is injective. Therefore, given a compact orientable relatively spin Lagrangian
submanifold $L\subset (M,\omega)$ such that $i_*:H_*(L,\Q)\to H_*(M,\Q)$ is
injective, for any Hamiltonian diffeomorphism $\phi:(M,\omega)\to
(M,\omega)$ such that $\phi(L)$ is transverse to $L$, we have 
\begin{equation}\#(L\cap \phi(L))\ge \sum_k \mathrm{rank}\, H_k(L,\Q).
\end{equation}
For example, consider
a symplectic manifold of the form $(M\times M,\omega\oplus -\omega)$, and
let $L_0$ be the diagonal, and $L_1$ be the graph of a Hamiltonian
diffeomorphism $\psi:(M,\omega)\to(M,\omega)$. Then $L_0$ and $L_1$ are
Lagrangian submanifolds, and (3) yields the Arnold conjecture for the fixed
points of $\psi$.

Besides a differential, Floer complexes for Lagrangians are also equipped
with a product structure, i.e.\ a morphism $CF^*(L_0,L_1)\otimes
CF^*(L_1,L_2)\to CF^*(L_0,L_2)$ (assuming that the obstruction classes
vanish for the given Lagrangians $L_0,L_1,L_2$).
This product structure is defined as follows:
consider three points $p_1\in L_0\cap L_1$, $p_2\in L_1\cap L_2$,
$p_3\in L_0\cap L_2$, together with the equivalence classes
$[u_i]$ needed to specify generators of $CF^*$. Consider the
moduli space $\mathcal{M}(p_1,p_2,p_3;A)$ of all pseudo-holomorphic
maps $f$ from a disc with three marked points $q_1,q_2,q_3$ on its
boundary to $M$, such that $f(q_i)=p_i$ and the three portions of boundary
delimited by the marked points are mapped to $L_0,L_1,L_2$ respectively,
and realizing a given relative homology class $A$. We compactify this moduli
space and complete it if necessary in order to obtain a well-defined
fundamental cycle. Assuming that the relative homology class $A$ is
compatible with the given equivalence classes $[u_i]$, the virtual
dimension of this moduli space is the difference between the Maslov index
of $(p_3,[u_3])$ and the sum of
those of $(p_1,[u_1])$ and $(p_2,[u_2])$ . The product of 
$(p_1,[u_1])$ and $(p_2,[u_2])$ is then defined as
\begin{equation}
(p_1,[u_1])\cdot (p_2,[u_2])=\sum_{p_3,A} \#\mathcal{M}(p_1,p_2,p_3;A)\,
(p_3,[u_3]),\end{equation}
where the sum runs over all pairs $(p_3,A)$ for which the expected dimension
of the moduli space is zero, and for each such pair $[u_3]$ is the
equivalence class determined by $[u_1]$, $[u_2]$ and $A$.

While the product structure on $CF^*$ defined by (4) satisfies the Leibniz
rule with respect to the differential $\partial$ (and hence descends to a
product structure on Floer homology), it differs from usual products
by the fact that it is only associative {\it up to homotopy}. In fact, Floer
complexes come equipped with a full set of {\it higher-order products}
$$\mu^n:CF^*(L_0,L_1)\otimes \cdots\otimes CF^*(L_{n-1},L_n)\to CF^*(L_0,L_n)
\quad \mbox{for all}\ n\ge 1,$$
with each $\mu^n$ shifting degree by $2-n$. The first two maps $\mu^1$ and
$\mu^2$ are respectively the Floer differential $\partial$ and the
product described above. The definition of $\mu^n$ is similar to
those of $\partial$ and of the product structure: given
generators $(p_i,[u_i])\in CF^*(L_{i-1},L_i)$ for $1\le i\le n$ and
$(p_{n+1},[u_{n+1}])\in CF^*(L_0,L_n)$ such that $\deg (p_{n+1},[u_{n+1}])=
\sum_{i=1}^n \deg (p_i,[u_i])+2-n$, the coefficient of $(p_{n+1},[u_{n+1}])$
in $\mu^n((p_1,[u_1]),\dots,(p_n,[u_n])$ is obtained by counting (in a
suitable sense) pseudo-holomorphic 
maps $f$ from a disc with $n+1$ marked points $q_1,\dots,q_{n+1}$ on its
boundary to $M$, such that $f(q_i)=p_i$ and the portions of boundary
delimited by the marked points are mapped to $L_0,\dots,L_n$ respectively,
and representing a relative homology class compatible with the
given data $[u_i]$.

Assume that there are no non-trivial pseudo-holomorphic discs with boundary
contained in one of the considered Lagrangian submanifolds (otherwise
one needs to introduce a $\mu^0$ term as well, which is perfectly legitimate
but makes the structure much more different from that of a usual category): then 
the maps $(\mu^n)_{n\ge 1}$ define an {\em $A_\infty$-structure} on Floer
complexes, i.e.\ they satisfy an infinite sequence of algebraic relations:
\smallskip
\begin{equation*}\begin{cases}
\mu^1(\mu^1(a))=0,\\
\mu^1(\mu^2(a,b))=\mu^2(\mu^1(a),b)+(-1)^{\deg a}\mu^2(a,\mu^1(b)),\\
\mu^1(\mu^3(a,b,c))=\mu^2(\mu^2(a,b),c)-\mu^2(a,\mu^2(b,c))\\
\qquad\qquad\qquad\qquad
\pm\mu^3(\mu^1(a),b,c)\pm\mu^3(a,\mu^1(b),c)\pm\mu^3(a,b,\mu^1(c)),\\
\cdots
\end{cases}\end{equation*}

This leads to the concept of ``Fukaya category'' of a symplectic manifold.
Conjecturally, for every symplectic manifold $(M,\omega)$ one should be able
to define an \hbox{$A_\infty$-category} $\mathcal{F}(M)$ whose objects are
Lagrangian submanifolds (compact, orientable, relatively spin, ``twisted''
by a flat unitary vector bundle);
the space of morphisms between two objects $L_0$ and $L_1$ is the Floer
complex $CF^*(L_0,L_1)$ equipped with its differential $\partial=\mu^1$, with
(non-associative) composition given by the product $\mu^2$, and higher order
compositions $\mu^n$.

The importance of Fukaya categories in modern symplectic topology is largely
due to the {\it homological mirror symmetry} conjecture, formulated by
Kontsevich. Very roughly, this conjecture states that the phenomenon of mirror
symmetry, i.e.\ a conjectural correspondence between symplectic manifolds
and complex manifolds (``mirror pairs'') arising from a duality among string
theories, should be visible at the level of
Fukaya categories of symplectic manifolds and categories of coherent sheaves
on complex manifolds: given a mirror pair consisting of a symplectic
manifold $M$ and a complex manifold $X$, the derived categories
$D\mathcal{F}(M)$ and $D^b Coh(X)$ should be equivalent (in a more precise form
of the conjecture, one should actually consider families of manifolds and
deformations of categories). However, due to the very
incomplete nature of our understanding of Fukaya categories in comparison to
the much better understood derived categories of coherent sheaves, this
conjecture has so far only been verified on very specific examples.

\subsection{The topology of symplectic 4-manifolds}

An important question in symplectic topology is to determine which smooth
manifolds admit symplectic structures. In the case of open manifolds,
Gromov's $h$-principle implies that the existence of an almost-complex
structure is sufficient. In contrast, the case of compact manifolds is much
less understood, except in dimension 4.

Whereas the existence of a class $\alpha\in H^2(M,\R)$ such that
$\alpha^{\cup n}\neq 0$ and of an almost-complex structure already provide
elementary obstructions to the existence of a symplectic structure on a
given compact manifold, in the case of 4-manifolds a much stronger
obstruction arises from {\it Seiberg-Witten invariants}.

The Seiberg-Witten equations are defined on a compact 4-manifold $M$ equipped
with a Riemannian metric and a $\mathrm{spin}^c$ structure $\mathfrak{s}$,
characterized by a pair of rank $2$ complex vector bundles $S^\pm$
(positive and negative spinors, interchanged by the {\it Clifford action} of
the tangent bundle), with the same determinant line bundle $L$.
The choice of a Hermitian connection $A$ on $L$ determines (together with
the Levi-Civita connection on the tangent bundle) a $\mathrm{spin}^c$
connection on $S^\pm$, which in turn yields a Dirac operator
$D_A:\Gamma(S^\pm)\to\Gamma(S^\mp)$ (obtained by contraction of the connection
operator with the Clifford action); moreover, the bundle of self-dual 2-forms
$\Lambda^2_+ T^*M$ is canonically isomorphic (via the Clifford action) to
that of traceless antihermitian endomorphisms of $S^+$. With this
understood, the Seiberg-Witten equations are the following equations for a
pair $(A,\psi)$
consisting of a $U(1)$ connection on the determinant line bundle $L$ and a
section of the bundle of positive spinors \cite{Wi}:
\begin{equation}
\begin{cases}
D_A\psi=0\\
F_A^+=q(\psi)+\mu
\end{cases}
\end{equation}
where $q(\psi)$ is the imaginary self-dual 2-form corresponding to the
traceless part of $\psi^*\otimes \psi\in \mathrm{End}(S^+)$ and $\mu$ is
a constant parameter. The space of solutions to (5) is invariant under 
$U(1)$-valued gauge transformations $(A,\psi)\mapsto (A+2g^{-1}dg, g\psi)$,
and the quotient moduli space $\mathcal{M}(L)$ is compact,
orientable, smooth for generic choice of the parameter $\mu$ whenever
$b_2^+(M)\neq 0$, and of expected dimension
$d=\frac{1}{4}(c_1(L)^2-2\chi(M)-3\sigma(M))$; the Seiberg-Witten invariants
of $X$ are then defined by counting points of the moduli space (with signs):
when $d=0$ (the most important case) we set $SW(L)=\#\mathcal{M}(L)$.
Whenever $b_2^+(M)\ge 2$ this invariant depends only on the manifold $M$ and
the given $\mathrm{spin}^c$-structure; when $b_2^+(M)=1$ there
are two different chambers depending on the choice of the metric and the
parameter $\mu$, giving rise to possibly different values of the invariant.

The most important results concerning the Seiberg-Witten invariants of
symplectic 4-manifolds have been obtained by Taubes (\cite{Ta1},
\cite{Ta2}, \dots). They are summarized in the following statement:

\begin{theorem}[Taubes]
Let $(M^4,\omega)$ be a compact symplectic 4-manifold with $b_2^+\ge 2$.
Then:

$(i)$ $SW(K_M^{\pm 1})=\pm 1$;

$(ii)$ $c_1(K_M)\cdot [\omega]\ge 0$, and $SW(L)=0$
whenever $|c_1(L)\cdot [\omega]|>c_1(K_M)\cdot [\omega]$.

$(iii)$ $SW(K_M^{-1}+2e)=Gr_T(e)$, where $Gr_T$ is a specific version of
Gromov-Witten invariants (counting possibly disconnected pseudo-holomorphic
curves, with special weights attributed to multiply covered square zero
tori);

$(iv)$ the homology class $c_1(K_M)$ admits a (possibly disconnected)
pseudo-holo\-mor\-phic representative, every component of which satisfies $g=1+[C_i]\cdot
[C_i]$. Hence, if $M$ is minimal i.e.\ contains no $(-1)$-spheres,
then $c_1(K_M)^2=2\chi(M)+3\sigma(M)\ge 0$.
\end{theorem}

These criteria prevent many 4-manifolds from admitting symplectic
structures, e.g.\ those which decompose as the connected sum of two
manifolds with $b_2^+\ge 1$ (by a result of Witten, their SW invariants
vanish): for example, $\CP^2\#\CP^2\#\CP^2$ does not admit any symplectic
structure even though it satisfies the cohomological condition for the
existence of an almost-complex structure (there exists a class $c\in
H^2(M,\Z)$ such that $c^2=2\chi+3\sigma=19$).

When $b_2^+(M)=1$, some of the statements remain valid (being careful about
the choice of chamber for SW invariants). Using Gromov's characterization of
the Fubini-Study symplectic structure of $\CP^2$ in terms of the existence of
pseudo-holomorphic lines, Taubes has shown that the symplectic structure of
$\CP^2$ is unique up to scaling. This result has been extended by Lalonde
and McDuff to the case of rational ruled surfaces, where $\omega$ is
determined by its cohomology class.
\medskip

In parallel to the above constraints on symplectic 4-manifolds, surgery
techniques have led to many interesting new examples of compact symplectic
manifolds.

One of the most efficient techniques in this respect is the {\it symplectic
sum} construction, investigated by Gompf \cite{Go1}:
if two symplectic manifolds $(M_1^{2n},\omega_1)$ and $(M_2^{2n},\omega_2)$
contain compact symplectic hypersurfaces $W_1^{2n-2},W_2^{2n-2}$ that are
mutually symplectomorphic and whose normal bundles have opposite Chern
classes, then we can cut $M_1$ and $M_2$ open along the submanifolds $W_1$
and $W_2$, and glue them to each other along their common boundary,
performing a fiberwise connected sum in the normal bundles to $W_1$ and
$W_2$, to obtain a new symplectic manifold $M=M_1{\vphantom{M}}_{W_1}{\#}
{\vphantom{M}}_{W_2} M_2$. This construction has in particular allowed Gompf
to show that every finitely presented group can be realized as the
fundamental group of a compact symplectic 4-manifold.

A more specifically 4-dimensional construction is the {\it link surgery}
technique developed by Fintushel and Stern \cite{FS1}. By modifying a given
4-manifold in the neighborhood of an embedded torus with trivial normal
bundle, it allows one to construct large families of mutually homeomorphic
4-manifolds whose diffeomorphism types can be distinguished using
Seiberg-Witten invariants; these manifolds are symplectic in some cases.
More can be found about 4-dimensional surgery techniques in the book
\cite{GS}; we will also discuss some symplectic surgery constructions
in more detail in the sections below.

By comparing the available examples and the topological constraints arising
from Seiberg-Witten invariants, one can get a fairly good understanding of
the topology of compact symplectic 4-manifolds. Nonetheless, a large
number of questions remain open, in particular concerning symplectic
manifolds of {\it general type}, i.e.\ with $b_2^+\ge 2$ and $c_1(K_M)^2>0$.
For example, Seiberg-Witten invariants fail to provide any useful
information on complex surfaces of general type, whose diffeomorphism types
are hence not well understood; similarly, it is unknown to this date whether
the Bogomolov-Miyaoka-Yau inequality $c_1^2\le 3 c_2$, satisfied by all
complex surfaces of general type, also holds in the symplectic case.

%}

\section{Symplectic Lefschetz fibrations}

%\comment{

\subsection{Symplectic fibrations}

The first construction of a symplectic non-K\"ahler manifold is due to
Thurston \cite{Th}, who showed the existence of a symplectic structure on
a $T^2$-bundle over $T^2$ with $b_1=3$ (whereas K\"ahler manifolds always
have even $b_1$ since their first cohomology splits into $H^{1,0}\oplus
H^{0,1}$). The existence of a symplectic structure on this manifold follows
from a general result about symplectic fibrations (see also \S 6 of
\cite{McS}):

\begin{theorem}[Thurston]
Let $f:M\to B$ be a compact locally trivial fibration with symplectic
fiber $(F,\omega_F)$ and symplectic base $(B,\omega_B)$. Assume that
the structure group of the fibration reduces to the symplectomorphisms of
$F$, and that there exists a cohomology class $c\in H^2(M,\R)$ whose
restriction to the fiber is equal to $[\omega_F]$. Then, for all sufficiently
large $K>0$, $M$ admits a symplectic form in the cohomology class
$c+K\,f^*[\omega_B]$, for which all fibers of $f$ are symplectic
submanifolds.
\end{theorem}

When the fibers are $2$-dimensional, many of the assumptions of this theorem
are always satisfied: the fiber always admits a symplectic structure
(provided it is orientable), and the structure group of the fibration
always reduces to symplectomorphisms; moreover, the cohomological condition
is equivalent to the requirement that the fibers of $f$ represent a non-zero
class in $H_2(M,\R)$.

\proof
Let $\eta\in\Omega^2(M)$ be a closed 2-form representing the cohomology
class $c$. Cover the base $B$ by balls $U_i$ over which the fibration
is trivial: we have a diffeomorphism $\phi_i:f^{-1}(U_i)\to
U_i\times F$, and the assumption on the structure group of the
fibration means that, over $U_i\cap U_j$, the trivializations
$\phi_i$ and $\phi_j$ differ by symplectomorphisms of the fibers.
The diffeomorphism $\phi_i$ determines a projection $p_i:f^{-1}(U_i)\to F$
such that $\phi_i(x)=(f(x),p_i(x))$ for every $x\in f^{-1}(U_i)$.

After restriction to $f^{-1}(U_i)\simeq U_i\times F$, the closed 2-forms
$\eta$ and $p_i^*\omega_F$ belong to the same cohomology class, so that we
can write $p_i^*\omega_F=\eta+d\alpha_i$ for some 1-form $\alpha_i$
over $f^{-1}(U_i)$.
Let $\{\rho_i\}$ be a partition of unity by smooth functions
$\rho_i:B\to [0,1]$ supported over the balls $U_i$, and let
$\tilde\eta=\eta+\sum_i d((\rho_i\circ f)\,\alpha_i)$.
The 2-form $\tilde\eta$ is well-defined since the support of $\rho_i\circ f$
is contained in $f^{-1}(U_i)$, and it is obviously closed and cohomologous
to $\eta$. Moreover, over $F_p=f^{-1}(p)$ for any $p\in B$, we have
$\tilde\eta_{|F_p}=\eta_{|F_p}+\sum_i \rho_i(p)
d\alpha_{i|F_p}=\sum_i \rho_i(p)\,(\eta+d\alpha_i)_{|F_p}=
\sum_i \rho_i(p)\,(p_i^*\omega_F)_{|F_p}$. However since the 
local trivializations $\phi_i$ differ by symplectomorphisms of the fiber,
the 2-forms $p_i^*\omega_F$ are all equal to each other. Therefore,
the 2-form $\tilde\eta_{|F_p}$ can be identified with $\omega_F$ (recall
that $\sum \rho_i\equiv 1$). 

So far we have constructed a closed 2-form $\tilde\eta\in
\Omega^2(M)$, with $[\tilde\eta]=c$, whose restriction to any fiber of $f$ is
symplectic (and in fact coincides with $\omega_F$). At any point $x\in M$,
the tangent space $T_xM$ splits into a vertical subspace
$V_x=\mathrm{Ker}\,df_x$ and a horizontal subspace $H_x=\{v\in T_xM,\ 
\tilde\eta(v,v')=0\ \forall v'\in V_x\}$. Since the restriction of
$\tilde\eta$ to the vertical subspace is non-degenerate, we have
$T_xM=H_x\oplus V_x$, and so $f^*\omega_B$ is non-degenerate over
$H_x$. Therefore, for sufficiently large $K>0$ the 2-form
$\tilde\eta+K\,f^*\omega_B$ is non-degenerate over $H_x$, and since its
restriction to $V_x$ coincides with $\tilde\eta$, it is also non-degenerate
over $T_xM$. Using the compactness of $M$ we can find a constant $K>0$ for
which this property holds at every point; the form
$\tilde\eta+K\,f^*\omega_B$ then defines a symplectic structure on $M$.
\endproof

Thurston's example of a non-K\"ahler symplectic manifold is obtained in the
following way: start with the trivial bundle with fiber $T^2=\R^2/\Z^2$
(with coordinates $x,y\in\R/\Z$) over $\R^2$ (with coordinates $z,t$), and
quotient it by the action of $\Z^2$ generated by $((x,y),(z,t))\mapsto
((x,y),(z+1,t))$ and $((x,y),(z,t))\mapsto ((x+y,y),(z,t+1))$. This action
of $\Z^2$ maps fibers to fibers, and hence the quotient $M$ carries a
structure of $T^2$-fibration over $\R^2/\Z^2=T^2$. The fiber class is
homologically non-trivial since it has intersection number $1$ with the
section $\{x=y=0\}$, and the monodromy of the fibration is given by
symplectomorphisms of the fiber ($(x,y)\mapsto (x,y)$ and $(x,y)\mapsto
(x+y,y)$). Therefore, by Theorem 2.1 the compact 4-manifold $M$ admits a
symplectic structure. Let $\gamma_1,\dots,\gamma_4$ be the closed loops in
$M$ corresponding to the four coordinate axes: by translating $\gamma_2$
along the $t$ axis, one can deform it into a loop homologous to
$\gamma_1+\gamma_2$ in the fiber, so that $[\gamma_1]=0$ in $H_1(M,\Z)$.
On the other hand, it is not hard to see that
$[\gamma_2],[\gamma_3],[\gamma_4]$ are linearly independent and generate
$H_1(M,\Z)\simeq \Z^3$. Hence $b_1(M)=3$ and $M$ does not admit any K\"ahler
structure (on the other hand, $M$ actually carries an integrable complex
structure, but it is not compatible with any symplectic form).

\subsection{Symplectic Lefschetz fibrations}

\begin{definition}
A map $f$ from a compact oriented manifold $M^{2n}$ to the sphere $S^2$ (or
more generally a compact oriented Riemann surface) is a
{\em Lefschetz fibration} if the critical points of $f$ are isolated,
and for every critical point $p\in M$ the map $f$ is modelled on a complex
Morse function, i.e.\ there exist neighborhoods $U\ni p$ and $V\ni f(p)$ and
orientation-preserving local diffeomorphisms
$\phi:U\to\C^n$ and $\psi:V\to\C$ such that $\psi\circ f\circ \phi^{-1}$ is
the map $(z_1,\dots,z_n)\mapsto \sum z_i^2$.
\end{definition}

To simplify the description, we can additionally require that the critical
values of $f$ are all distinct (so that each fiber contains at most one
singular point).

The local model near a singular fiber is easiest to understand in the
4-dimensional case. The fiber $F_\lambda$ of the map
$(z_1,z_2)\to z_1^2+z_2^2$ above $\lambda\in\C$ is given by the equation
$(z_1+iz_2)(z_1-iz_2)=\lambda$: the fiber $F_\lambda$ is smooth 
(topologically a cylinder) for all $\lambda\neq 0$, while the fiber
above the origin presents a
transverse double point, and is obtained from the nearby fibers by
collapsing an embedded simple closed loop called the {\it vanishing cycle}.
For example, for $\lambda>0$ the vanishing cycle is the loop
$\{(x_1,x_2)\in\R^2,\ x_1^2+x_2^2=\lambda\}=F_\lambda\cap \R^2\subset
F_\lambda$. In arbitrary dimension, the fiber over a critical point of $f$
presents an {\it ordinary double point}, and the nearby fibers are
smoothings of this singularity; this can be seen in the local model, where
the singular fiber $F_0=\{\sum z_i^2=0\}\subset \C^n$ is obtained from a
nearby smooth fiber $F_\lambda=\{\sum z_i^2=\lambda\}$ by
collapsing the vanishing cycle $S_\lambda=F_\lambda\cap
(e^{i\theta/2}\R)^n\subset F_\lambda$ (where $\theta=\arg(\lambda)$).
In fact, $S_\lambda$ is obtained from the unit sphere in $\R^n\subset\C^n$
by multiplication by $\lambda^{1/2}$, and $F_\lambda$ is
diffeomorphic to the cotangent bundle $T^*S_\lambda$.

Fix a base point $q_0\in S^2-\mathrm{crit}(f)$, and consider a closed loop
$\gamma:[0,1]\to S^2-\mathrm{crit}(f)$ (starting and ending at $q_0$). 
By fixing a horizontal distribution we can perform parallel transport in
the fibers of $f$ along $\gamma$, which induces a diffeomorphism from
$F_{q_0}=f^{-1}(q_0)$ to itself. The isotopy class of this diffeomorphism, which
is well-defined independently of the chosen horizontal
distribution, is called the {\it monodromy} of $f$ along $\gamma$. Hence, we
obtain a monodromy homomorphism characteristic of the Lefschetz fibration $f$,
$$\psi:\pi_1(S^2-\mathrm{crit}(f),q_0)\to
\pi_0\mathrm{Diff}^+(F_{q_0}).$$

By considering the local model near a critical point of $f$, one can show
that the monodromy of $f$ around one of its singular fibers is a
{\it positive Dehn twist} along the vanishing cycle -- a diffeomorphism
supported in a neighborhood of the vanishing cycle, and inducing
the antipodal map on the vanishing cycle itself (recall from above
that in the local model we have $S_\lambda=F_\lambda\cap
(e^{i\theta/2}\R)^n$). The effect of this Dehn twist is most easily
seen when $\dim M=4$: in the fiber $F_{q_0}$, a tubular neighborhood of the
vanishing cycle can be identified with a cylinder $\R\times S^1$, and
the opposite ends of the cylinder twist relatively to each other as one
moves around the singular fiber, as shown in the figure below. 

\begin{center}
\epsfig{file=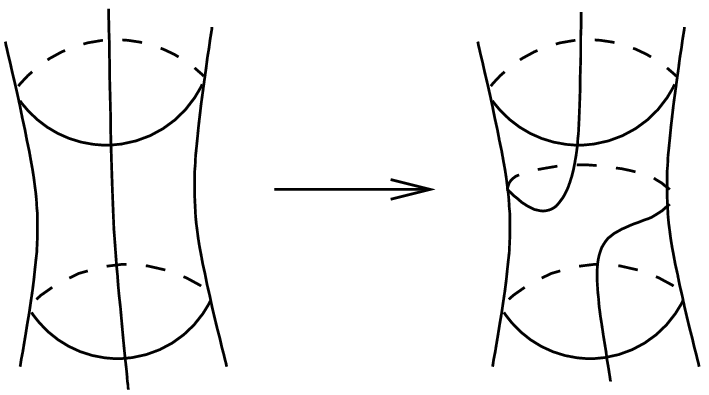,height=2cm,width=4.5cm}
\end{center}

In higher dimension, a neighborhood of the vanishing cycle is diffeomorphic
to the neighborhood $\{(x,u)\in T^*S^{n-1},\ |u|\le \pi\}$ of the zero
section in the cotangent bundle of $S^{n-1}$, and the Dehn twist is the
diffeomorphism supported in this neighborhood obtained by composing the
time $1$ geodesic flow on $T^*S^{n-1}$ with the differential of the
antipodal involution of $S^{n-1}$.

We are particularly interested in the case where the total space $M$ carries
a symplectic structure compatible with the fibration:

\begin{definition} A {\em symplectic Lefschetz fibration} $f:(M,\omega)\to
S^2$ is a Lefschetz fibration whose total space carries a symplectic
structure with the following properties: $(i)$ the fibers of $f$ are
symplectic submanifolds; $(ii)$ near a critical point of $f$, there exist
neighborhoods $U\ni p$ and $V\ni f(p)$ and orientation-preserving local
diffeomorphisms $\phi:U\to\C^n$ and $\psi:V\to\C$, such that the symplectic
form $\phi_*\omega$ evaluates positively on all complex lines in $\C^n$ 
$(\phi_*\omega(v,iv)>0$ $\forall v\neq 0)$, and such that
$\psi\circ f\circ \phi^{-1}$ is the map $(z_1,\dots,z_n)\mapsto \sum z_i^2$.
\end{definition}

By considering at each point
$p\in M$ where $df\neq 0$ the subspace of $T_pM$ symplectically orthogonal
to the fiber through $p$, one obtains a
specific horizontal distribution for which parallel transport preserves
the symplectic structures of the fibers of $f$. Hence, the monodromy of a
{\it symplectic Lefschetz fibration} can be defined with values into the
{\it symplectic mapping class group}, i.e.\ the group
$\pi_0\mathrm{Symp}(F)$ of isotopy classes of symplectomorphisms of the
fiber. Moreover, the vanishing cycles are always {\it Lagrangian} spheres
in the fibers (this is e.g.\ clearly true in the local model in $\C^n$
equipped with its standard symplectic structure; more generally, it follows
easily from the fact that parallel transport preserves the symplectic
structure and collapses the vanishing cycle into the critical point of a
singular fiber).

In the case of Lefschetz fibrations over a disc, the monodromy homomorphism
is sufficient to reconstruct the total space of the fibration up to
diffeomorphism (resp.\ symplectic deformation). When
considering fibrations over $S^2$, the monodromy data determines the
fibration over a large disc $D$ containing all critical values, after
which we only need to add a trivial fibration over a small disc $D'=S^2-D$,
to be glued in a manner compatible with the fibration structure over the
common boundary $\partial D=\partial D'=S^1$. This gluing involves the
choice of a map from $S^1$ to the group of diffeomorphisms
(resp.\ symplectomorphisms) of the fiber, i.e.\ an element of
$\pi_1\mathrm{Diff}(F)$ or $\pi_1\mathrm{Symp}(F)$.

In particular, in the 4-dimensional case, the total space of a (symplectic)
Lefschetz fibration with fibers of genus $g\ge 2$ is completely determined
by its monodromy, up to diffeomorphism (resp.\ symplectic deformation).

In view of the above discussion, it is not very surprising that a result
similar to Theorem 2.1 holds for Lefschetz fibrations \cite{GS,Go2}:

\begin{theorem}[Gompf]
Let $f:M^{2n}\to S^2$ be a Lefschetz fibration with symplectic
fiber $(F,\omega_F)$. Assume that the structure group of the fibration
reduces to the symplectomorphisms of $F$, and in particular that the
vanishing cycles are embedded Lagrangian spheres in $(F,\omega_F)$.
Assume moreover that there exists a cohomology class $c\in H^2(M,\R)$ whose
restriction to the fiber is equal to $[\omega_F]$; in the case $n=2$, the
cohomology class $c$ is also required to evaluate positively on every
component of every reducible singular fiber. Then, for all sufficiently
large $K>0$, $M$ admits a symplectic form in the cohomology class
$c+K\,f^*[\omega_B]$, for which all fibers of $f$ are symplectic
submanifolds.
\end{theorem}

The proof of Theorem 2.4 is essentially a refinement of the argument given for
Theorem 2.1. The main difference is that one first uses the local models near
the critical points of $f$ in order to build an exact perturbation
of $\eta$ with the desired behavior near these points; in the subsequent
steps, one then works with differential forms with support contained in the
complement of small balls centered at the critical points \cite{Go2}.

An important motivation for the study of symplectic Lefschetz fibrations is
the fact that, up to blow-ups, every compact symplectic manifold carries
such a structure, as shown by Donaldson \cite{Do2,Do3}:

\begin{theorem}[Donaldson]
For any compact symplectic manifold $(M^{2n},\omega)$, there exists a smooth
codimension $4$ symplectic submanifold $B\subset M$ such that the blow-up of
$M$ along $B$ carries a structure of symplectic Lefschetz fibration
over $S^2$.
\end{theorem}

In particular, this result (together with the converse statement of Gompf)
provides a very elegant topological description of symplectic 4-manifolds.
In order to describe Donaldson's construction of symplectic Lefschetz
pencils, we need a digression into
{\it approximately holomorphic geometry}.

\subsection{Approximately holomorphic geometry}

The idea introduced by Donaldson in the mid-90's is the following: in
presence of an almost-complex sructure, the lack of integrability usually
prevents the existence of holomorphic sections of vector bundles or
pseudo-holomorphic maps to other manifolds, but one can work in a similar manner
with approximately holomorphic objects.

Let $(M^{2n},\omega)$ be a compact symplectic manifold of dimension $2n$. We
will assume throughout this paragraph that $\frac{1}{2\pi}[\omega]\in
H^2(M,\Z)$; this integrality condition does not restrict the topological
type of $M$, since any symplectic form can be perturbed into another
symplectic form $\omega'$ whose cohomology class is rational (we can then
achieve integrality by multiplication by a constant factor). Morever, it is
easy to check that the submanifolds of $M$ that we will construct are not
only $\omega'$-symplectic but also $\omega$-symplectic, hence making the
general case of Theorem 2.5 follow from the integral case.

Let $J$ be an almost-complex structure compatible with $\omega$, and
let $g(.,.)=\omega(.,J.)$ be the corresponding Riemannian metric.
We consider a complex line bundle $L$ over $M$ such that
$c_1(L)=\frac{1}{2\pi}[\omega]$, endowed with a Hermitian metric and
a Hermitian connection $\nabla^L$ with curvature 2-form 
$F(\nabla^L)=-i\omega$. The almost-complex structure induces a
splitting of the connection~:
$\nabla^L=\partial^L+\dbar^L$, where $\partial^L s(v)=\frac{1}{2}(\nabla^L
s(v)-i\nabla^L s(Jv))$ and $\dbar^L s(v)=\frac{1}{2}(\nabla^L
s(v)+i\nabla^L s(Jv))$.

If the almost-complex structure $J$ is integrable, i.e.\ if $M$ is a 
K\"ahler complex manifold, then $L$ is an ample holomorphic line bundle,
and for large enough values of $k$ the line bundles $L^{\otimes k}$ admit
many holomorphic sections. Therefore, the manifold $M$ can be embedded into
a projective space (Kodaira); generic hyperplane sections are smooth
hypersurfaces in $M$ (Bertini), and more generally the linear system formed
by the sections of $L^{\otimes k}$ allows one to construct various
structures (Lefschetz pencils, \dots\!\!).

When the manifold $M$ is only symplectic, the lack of integrability of $J$
prevents the existence of holomorphic sections. Nonetheless, it is possible
to find an {\it approximately holomorphic} local model: a neighborhood of a
point $x\in M$, equipped with the symplectic form $\omega$ and the
almost-complex structure $J$, can be identified with a neighborhood of the
origin in $\C^n$ equipped with the standard symplectic form $\omega_0$ and
an almost-complex structure of the form $i+O(|z|)$. In this local model, the
line bundle $L^{\otimes k}$ endowed with the connection
$\nabla=(\nabla^L)^{\otimes k}$ of curvature $-ik\omega$ can be identified
with the trivial line bundle $\underline{\C}$ endowed with the connection 
$d+\frac{k}{4}\sum (z_j\,d\bar{z}_j-\bar{z}_j\,dz_j)$. The section of 
$L^{\otimes k}$ given in this trivialization by $s_{k,x}(z)=
\exp(-\frac{1}{4}k|z|^2)$ is then approximately holomorphic \cite{Do1}.

More precisely, a sequence of sections $s_k$ of $L^{\otimes k}$ is
said to be approximately holomorphic if, with respect to the rescaled
metrics $g_k=kg$, and after normalization of the sections to ensure that
$\|s_k\|_{C^r,g_k}\sim C$, an inequality of the form
$\|\dbar s_k\|_{C^{r-1},g_k}<C'k^{-1/2}$ holds, where $C$ and $C'$ are
constants independent of $k$. The change of metric, which dilates all
distances by a factor of $\smash{\sqrt{k}}$, is required in order to be able to
obtain uniform estimates, due to the larger and larger curvature of
the line bundle $L^{\otimes k}$. The intuitive idea is that, for large $k$,
the sections of the line bundle $L^{\otimes k}$ with curvature
$-ik\omega$ probe the geometry of $M$ at small scale
($\sim 1/\sqrt{k}$), which makes the almost-complex structure
$J$ almost integrable and allows one to achieve better and better
approximations of the holomorphicity condition $\dbar s=0$.

It is worth noting that, since the above requirement is an open condition,
it is not possible to define a ``space of approximately holomorphic
sections'' of $L^{\otimes k}$ in any simple manner (cf.\ the work of
Borthwick and Uribe \cite{BU} or Shiffman and Zelditch for other approaches
to this problem).

Once many approximately holomorphic sections have been made available,
the aim is to find among them some sections whose geometric behavior is
as generic as possible. Donaldson has obtained the following result 
\cite{Do1}:

\begin{theorem}[Donaldson]
For $k\gg 0$, $L^{\otimes k}$ admits approximately holomorphic sections
$s_k$ whose zero sets $W_k$ are smooth symplectic hypersurfaces.
\end{theorem}

This result starts from the observation that, if the section $s_k$ vanishes
transversely and if $|\dbar s_k(x)|\ll |\partial s_k(x)|$ at every point of
$W_k=s_k^{-1}(0)$, then the submanifold $W_k$ is symplectic,
and even approximately $J$-holomorphic
(i.e.\ $J(TW_k)$ is close to $TW_k$). The crucial point is therefore to
obtain a lower bound for $\partial s_k$ at every point of $W_k$, in order
to make up for the lack of holomorphicity.

Sections $s_k$ of $L^{\otimes k}$ are said to be {\it
uniformly transverse to $0$} if there exists a constant $\eta>0$
(independent of $k$) such that the inequality $|\partial s_k(x)|_{g_k}>\eta$
holds at any point of $M$ where $|s_k(x)|<\eta$. In order to prove Theorem
2.6, it is sufficient to achieve this uniform estimate on the tranversality of
some approximately holomorphic sections $s_k$. The idea of the construction
of such sections consists of two main steps. The first one is an effective
local transversality result for complex-valued functions. Donaldson's
argument makes use of a result of Yomdin on the complexity of real
semi-algebraic sets; however a somewhat simpler argument can be used
instead \cite{Au5}. The second step is a remarkable globalization process,
which makes it possible to achieve uniform transversality over larger and
larger open subsets by means of successive perturbations of the sections
$s_k$, until transversality holds over the entire manifold $M$ \cite{Do1}.

The symplectic submanifolds constructed by Donaldson present several
remarkable properties which make them closer to complex submanifolds than
to arbitrary symplectic submanifolds. For instance, they satisfy the
Lefschetz hyperplane theorem: up to half the dimension of the submanifold,
the homology and homotopy groups of $W_k$ are identical to those of $M$
\cite{Do1}. More importantly, these submanifolds are, in a sense,
asymptotically unique: for given large enough $k$, the submanifolds 
$W_k$ are, up to symplectic isotopy, independent of all the choices made
in the construction (including that of the almost-complex structure $J$)
\cite{Au1}.

It is worth mentioning that analogues of Donaldson's
construction for contact manifolds have been obtained by
Ibort, Martinez-Torres and Presas (\cite{IMP}, \dots\!\!); see also recent
work of Giroux and Mohsen \cite{GM}.

\subsection{Symplectic Lefschetz pencils}

We now move on to Donaldson's construction of symplectic Lefschetz pencils
\cite{Do2,Do3}. In comparison with Theorem 2.6, the general setup is the
same, the main difference being that we consider no longer one, but
two sections of $L^{\otimes k}$. A pair of suitably chosen approximately
holomorphic sections $(s_k^0,s_k^1)$ of $L^{\otimes k}$ defines a family of
symplectic hypersurfaces $$\Sigma_{k,\alpha}=\{x\in M,\ s_k^0(x)-\alpha
s_k^1(x)=0\},\ \ \alpha\in\CP^1=\C\cup\{\infty\}.$$ The submanifolds 
$\Sigma_{k,\alpha}$ are all smooth except for finitely many of them
which present an isolated singularity; they intersect transversely along
the {\it base points} of the pencil, which form a smooth symplectic
submanifold $Z_k=$\hbox{$\{s_k^0=s_k^1=0\}$} of codimension $4$.

The two sections $s_k^0$ and $s_k^1$ determine a projective map
$f_k=(s_k^0:s_k^1):M-Z_k\to\CP^1$, whose critical points correspond to the
singularities of the fibers $\Sigma_{k,\alpha}$. In the case of a symplectic
Lefschetz pencil, the function $f_k$ is a
complex Morse function, i.e.\ near any of its critical points it is given
by the local model $f_k(z)=z_1^2+\dots+z_n^2$ in approximately holomorphic
coordinates. After blowing up $M$ along $Z_k$, the Lefschetz pencil
structure on $M$ gives rise to a well-defined map $\hat{f}_k:\hat{M}\to
\CP^1$; this map is a symplectic Lefschetz fibration.
Hence, Theorem 2.5 may be reformulated more precisely as follows:

\begin{theorem}[Donaldson]
For large enough $k$, the given manifold $(M^{2n},\omega)$ admits symplectic
Lefschetz pencil structures determined by pairs of suitably chosen
approximately holomorphic sections $s_k^0,s_k^1$ of $L^{\otimes k}$.
Moreover, for large enough $k$ these Lefschetz pencil structures are
uniquely determined up to isotopy.
\end{theorem}

As in the case of submanifolds, Donaldson's argument relies on successive
perturbations of given approximately holomorphic sections $s_k^0$ and
$s_k^1$ in order to achieve uniform transversality properties, not only for
the sections $(s_k^0,s_k^1)$ themselves but also for the derivative
$\partial f_k$ \cite{Do3}.

The precise meaning of the uniqueness statement is the following:
assume we are given two sequences of Lefschetz pencil structures on $(M,\omega)$,
determined by pairs of approximately holomorphic sections of $L^{\otimes k}$
satisfying uniform transversality estimates, but possibly with respect to
two different $\omega$-compatible almost-complex structures on $M$. Then,
beyond a certain (non-explicit) value of $k$, it becomes possible to find
one-parameter families of Lefschetz pencil structures interpolating between
the given ones. In particular, this implies that for large $k$ the monodromy
invariants associated to these Lefschetz pencils only depend on
$(M,\omega,k)$ and not on the choices made in the construction.

The monodromy invariants associated to a symplectic Lefschetz pencil are
essentially those of the symplectic Lefschetz fibration obtained after
blow-up along the base points, with only a small refinement. After the
blow-up operation, each fiber of $\hat{f_k}:\hat{M}\to \CP^1$ contains a
copy of the base locus $Z_k$ embedded as a smooth symplectic hypersurface.
This hypersurface lies away from all vanishing cycles, and is preserved by
the monodromy. Hence, the monodromy homomorphism can be defined to take
values in the group of isotopy classes of symplectomorphisms of the fiber
$\Sigma_k$ whose restriction to the submanifold $Z_k$ is the identity.

We can even do better by carefully examining the local model near the base
points and observing that parallel transport (with respect to the natural
symplectic connection) along a loop $\gamma\subset
\CP^1-\mathrm{crit}(\hat{f}_k)$ acts on a neighborhood of $Z_k$ by
complex multiplication by $e^{-2i\pi A(\gamma)}$ in the fibers of the normal
bundle, where $A(\gamma)$ is the proportion of the area of $\CP^1$
enclosed by $\gamma$. If we remove a regular value of $\hat{f}_k$
(e.g.\ the point at infinity in $\CP^1$), the area enclosed
by a given oriented loop becomes a well-defined element of $\R$ (rather
than $\R/\Z$), so that we can unambiguously represent every homotopy
class by a loop for which the enclosed area is equal to zero.
This makes it possible to define a monodromy homomorphism
\begin{equation}
\psi_k:\pi_1(\C-\mathrm{crit}(\hat{f}_k))\to \mathrm{Map}^\omega(\Sigma_k,
Z_k)\end{equation}
with values in the {\it relative} symplectic mapping class group
$$\mathrm{Map}^\omega(\Sigma_k,Z_k)=\pi_0(\{\phi\in \mathrm{Symp}(\Sigma_k,
\omega_{|\Sigma_k}),\ \phi_{|V(Z_k)}=\mathrm{Id}\}),$$ where
$\Sigma_k$ is a generic fiber of $\hat{f}_k$ (above the chosen
base point in $\C-\mathrm{crit}(\hat{f}_k)$), and $V(Z_k)$ is a neighborhood
of $Z_k$ inside $\Sigma_k$. As before, the monodromy around a singular fiber
is a positive Dehn twist along the vanishing cycle, which is
an embedded Lagrangian sphere $S^{n-1}\subset \Sigma_k$. However, the product
of all the monodromies around the individual singular fibers (i.e., the
``monodromy at infinity'') is not trivial, but rather equal to an element
$\delta_{Z_k}\in
\mathrm{Map}^\omega(\Sigma_k,Z_k)$, the positive Dehn twist
along the unit sphere bundle in the normal bundle of $Z_k$ in $\Sigma_k$
(i.e.\ $\delta_{Z_k}$ restricts to each fiber of the normal bundle as a Dehn
twist around the origin).

Things are easiest when $M$ is a 4-manifold: the fibers are then
Riemann surfaces, and $Z_k$ is a finite set of points.
The group $\mathrm{Map}^\omega(\Sigma_k,Z_k)$ therefore identifies with
the mapping class group $\mathrm{Map}_{g,N}$ of a Riemann surface of
genus $g=g(\Sigma_k)$ with $N=\mathrm{card}\,Z_k$ boundary components. The
monodromy around a singular fiber (a Riemann surface with an ordinary double
point) is a Dehn twist along an embedded loop, and the element $\delta_{Z_k}$
is the product of the Dehn twists along $N$ small loops encircling the
punctures.

Finally, we mention that Gompf's result (Theorem 2.4) admits an adaptation
to the case of Lefschetz pencils \cite{Go2}: in the presence of base points,
the cohomology class $[\omega]$
is determined in advance by the topology of the pencil, so that (using
Moser's theorem) the constructed symplectic form on blown-down manifold
$M$ becomes canonical
{\it up to symplectomorphism} (rather than just deformation equivalence).

%}

%\setcounter{section}{2}
%\setcounter{equation}{6}

\section{Lefschetz fibrations: examples and applications}

%\comment{

\subsection{Examples and classification questions} 
In this section, we give examples of Lefschetz fibrations, and mention some
results and questions related to their classification. We restrict ourselves
to the case where the total space is of dimension 4, because it is by far
the best understood. In fact, the structure of the symplectic mapping class
group $\mathrm{Map}^\omega(\Sigma)=\pi_0\mathrm{Symp}(\Sigma)$ of a symplectic
4-manifold is almost never known, except in the simplest cases (always
rational or ruled surfaces), which makes it very difficult to say much
about higher-dimensional symplectic Lefschetz fibrations (see \S 4 for an
alternative approach to this problem).

As a first example, we consider a pencil of
degree $2$ curves in $\CP^2$, a case in which the monodromy homomorphism (6)
can be directly explicited in a remarkably simple manner. 
Consider the two sections
$s_0=x_0(x_1-x_2)$ and $s_1=x_1(x_2-x_0)$ of the line bundle $O(2)$ over
$\CP^2$: their zero sets are singular conics, in fact the unions of two
lines each containing two of the four intersection points $(1\!:\!0\!:\!0)$,
$(0\!:\!1\!:\!0)$, $(0\!:\!0\!:\!1)$, $(1\!:\!1\!:\!1)$. Moreover, the
zero set of the linear combination $s_0+s_1=x_2(x_1-x_0)$ is also singular;
on the other hand, it is fairly easy to check that all other linear
combinations $s_0+\alpha s_1$ (for $\alpha\in \CP^1-\{0,1,\infty\}$) vanish
along smooth conics. By blowing up the four base points of the pencil
generated by $s_0$ and $s_1$, we
obtain a genus $0$ Lefschetz fibration with four exceptional sections.
The three singular fibers are nodal configurations consisting of
two tranversely intersecting spheres, with each component containing
two of the four base points; each of the three different
manners in which four points can be split into two groups of two is realized
at one of the singular fibers. The following diagram represents the three
singular conics of the pencil inside $\CP^2$ (left), and the corresponding
vanishing cycles inside a smooth fiber (right):
\medskip

\begin{center}
\epsfig{file=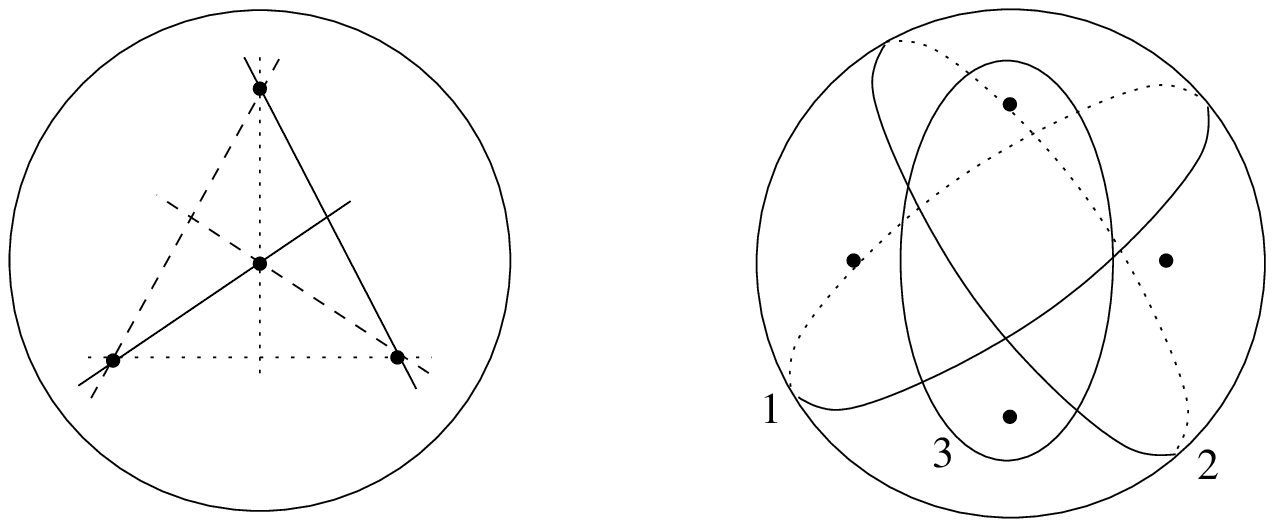,height=1.2in}
\end{center}
As seen above, the monodromy of this Lefschetz pencil can be expressed by
a morphism $\psi:\pi_1(\C-\{p_1,p_2,p_3\})\to \mathrm{Map}_{0,4}$. After
choosing
a suitable ordered basis of the free group $\pi_1(\C-\{p_1,p_2,p_3\})$,
we can make sure that $\psi$ maps the generators to the Dehn twists
$\tau_1,\tau_2,\tau_3$ along the three loops shown on the diagram.
On the other hand, as seen in \S 2.4 the monodromy at infinity is given
by the boundary twist $\prod \delta_i$, the product of the four Dehn twists
along small loops encircling the four base points in the fiber. The
monodromy at infinity can be decomposed into the product of the monodromies
around each of the three singular fibers ($\tau_1,\tau_2,\tau_3$). Hence, the
monodromy of a pencil of conics in $\CP^2$ can be expressed by the relation
$\prod \delta_i=\tau_1\cdot \tau_2\cdot \tau_3$ in the mapping class group
$\mathrm{Map}_{0,4}$ ({\it lantern relation}).

In our subsequent examples, we forget the base points,
and simply consider Lefschetz fibrations and their monodromy in
$\mathrm{Map}_g$ (rather than Lefschetz pencils with monodromy in
$\mathrm{Map}_{g,N}$).

The classification of genus $1$ Lefschetz fibrations over $S^2$ is a
classical result of Moishezon and Livne \cite{Mo1}, who have shown that,
up to isotopy, any such fibration is holomorphic. The mapping class group
$\mathrm{Map}_1$ is isomorphic to $SL(2,\Z)$ (the isomorphism is given by
considering the action of a diffeomorphism of $T^2$ on $H_1(T^2,\Z)=\Z^2$).
It is generated by the Dehn twists $$\tau_a=\left(\begin{array}{cc}1&1\\
0&1\end{array}\right)\quad\mathrm{and}\quad \tau_b=\left(\begin{array}{rc}
1&0\\\!\!-1&1\end{array}\right)$$
along the two generating loops $a=S^1\times\{pt\}$ and $b=\{pt\}\times S^1$,
with relations $\tau_a\tau_b\tau_a=\tau_b\tau_a\tau_b$ and
$(\tau_a\tau_b)^6=1$. If we consider a {\it relatively minimal} genus 1
Lefschetz fibration (i.e., if there are no homotopically trivial vanishing
cycles, an assumption that always becomes true after blowing down spherical
components of reducible singular fibers), then the number of singular fibers
is a multiple of 12, and for a suitable choice of an ordered system of
generating loops $\gamma_1,\dots,\gamma_{m=12k}$ in the base of the fibration,
the vanishing cycles can always be assumed to be $a,b,a,b,\dots$ In other
terms, the monodromy morphism maps the defining relation $\gamma_1\cdot\ldots
\cdot \gamma_m=1$ of the group $\pi_1(\CP^1-\mathrm{crit})$ to the positive
relation $(\tau_a\cdot\tau_b)^{6k}=1$ among Dehn twists in $\mathrm{Map}_1$.

The situation becomes more interesting for genus $2$ Lefschetz fibrations.
The mapping class group of a genus $2$ Riemann surface is generated by five
Dehn twists $\tau_i$ ($1\le i\le 5$) represented on the following diagram:
\medskip

\begin{center}
\setlength{\unitlength}{1cm}
\begin{picture}(5,2)(0,-1)
\qbezier[200](0,0)(0,1.2)(1.5,1)
\qbezier[200](1.5,1)(2.5,0.85)(3.5,1)
\qbezier[200](3.5,1)(5,1.2)(5,0)
\qbezier[200](0,0)(0,-1.2)(1.5,-1)
\qbezier[200](1.5,-1)(2.5,-0.85)(3.5,-1)
\qbezier[200](3.5,-1)(5,-1.2)(5,0)
\qbezier[60](1,0)(1,0.3)(1.5,0.3)
\qbezier[60](2,0)(2,0.3)(1.5,0.3)
\qbezier[60](1,0)(1,-0.3)(1.5,-0.3)
\qbezier[60](2,0)(2,-0.3)(1.5,-0.3)
\qbezier[60](3,0)(3,0.3)(3.5,0.3)
\qbezier[60](4,0)(4,0.3)(3.5,0.3)
\qbezier[60](3,0)(3,-0.3)(3.5,-0.3)
\qbezier[60](4,0)(4,-0.3)(3.5,-0.3)
\put(1.5,0){\circle{1.1}}
\put(3.5,0){\circle{1.1}}
\qbezier[60](0,0)(0,-0.1)(0.5,-0.1)
\qbezier[60](1,0)(1,-0.1)(0.5,-0.1)
\qbezier[20](0,0)(0,0.1)(0.5,0.1)
\qbezier[20](1,0)(1,0.1)(0.5,0.1)
\qbezier[60](2,0)(2,-0.1)(2.5,-0.1)
\qbezier[60](3,0)(3,-0.1)(2.5,-0.1)
\qbezier[20](2,0)(2,0.1)(2.5,0.1)
\qbezier[20](3,0)(3,0.1)(2.5,0.1)
\qbezier[60](4,0)(4,-0.1)(4.5,-0.1)
\qbezier[60](5,0)(5,-0.1)(4.5,-0.1)
\qbezier[20](4,0)(4,0.1)(4.5,0.1)
\qbezier[20](5,0)(5,0.1)(4.5,0.1)
%\qbezier[100](2.5,0.91)(2.6,0.91)(2.6,0)
%\qbezier[100](2.5,-0.91)(2.6,-0.91)(2.6,0)
%\qbezier[30](2.5,0.91)(2.4,0.91)(2.4,0)
%\qbezier[30](2.5,-0.91)(2.4,-0.91)(2.4,0)
\put(-0.35,0){$\tau_1$}
\put(0.8,0.5){$\tau_2$}
\put(2.1,-0.35){$\tau_3$}
\put(3.9,0.5){$\tau_4$}
\put(5.1,0){$\tau_5$}
%\put(2.65,0.7){$\sigma$}
\end{picture}
\end{center}

It is well-known (cf.\ e.g.\ \cite{Bi}, Theorem 4.8) that $\mathrm{Map}_2$
admits a presentation with generators $\tau_1,\dots,\tau_5$, and the
following relations: $\tau_i\tau_j=\tau_j\tau_i$ if $|i-j|\ge 2$;
$\tau_i\tau_{i+1}\tau_i=\tau_{i+1}\tau_i\tau_{i+1}$;
$(\tau_1\tau_2\tau_3\tau_4\tau_5)^6=1$; $I=\tau_1\tau_2\tau_3\tau_4\tau_5
\tau_5\tau_4\tau_3\tau_2\tau_1$ (the hyperelliptic involution) is central,
and $I^2=1$.

We now give two very classical examples of holomorphic genus $2$ Lefschetz
fibrations. The first one is the fibration $f_0:X_0\to\CP^1$, obtained from
a pencil of curves of bi-degree $(2,3)$ on $\CP^1\times\CP^1$ by blowing up
the 12 base points: this fibration has 20 singular fibers, and its monodromy
is described by the positive relation $(\tau_1\tau_2\tau_3\tau_4\tau_5\tau_5
\tau_4\tau_3\tau_2\tau_1)^2=1$ in $\mathrm{Map}_2$, expressing the identity
element as the product of the 20 Dehn twists along the vanishing cycles.
The second standard example of holomorphic genus 2 fibration is a fibration
$f_1:X_1\to\CP^1$ with total space a K3 surface blown up in two points, and
with 30 singular fibers; its monodromy is described by the positive relation
$(\tau_1\tau_2\tau_3\tau_4\tau_5)^6=1$ involving the Dehn twists along the
30 vanishing cycles.

If we have two Lefschetz fibrations $f:X\to \CP^1$ and $f':X'\to\CP^1$ with
fibers of the same genus, we can perform a {\it fiber sum} operation: choose
smooth fibers $F=f^{-1}(p)\subset X$ and $F'={f'}^{-1}(p')\subset X'$, and
an orientation-preserving diffeomorphism $\phi:F\to F'$. The complement
$U$ of a neighborhood of $F$ in $X$ is a symplectic 4-manifold with boundary
$F\times S^1$, and similarly the complement $U'$ of a neighborhood of $F'$ in
$X'$ has boundary $F'\times S^1$; by restriction to $U$ and $U'$, the
fibrations $f$ and $f'$ induce Lefschetz fibration structures over the
disc $D^2$. Using the diffeomorphism $\phi\times \mathrm{Id}$ to identify
$F\times S^1$ with $F'\times S^1$ in a manner compatible with the fibrations,
we can glue $U$ and $U'$ along their boundary in order to obtain
a compact symplectic manifold $X\#_\phi X'$ equipped with a Lefschetz
fibration $f\#_\phi f':X\#_\phi X'\to\CP^1$. In many cases there is a
particularly natural choice of gluing diffeomorphism $\phi$, leading to a
preferred fiber sum; fiber sums obtained using different gluing 
diffeomorphisms are then said to be {\it twisted}.

The fiber sum construction leads to many interesting examples of Lefschetz
fibrations; however, when the building blocks are the holomorphic
genus 2 fibrations $f_0$ and $f_1$, the result is actually independent of
the chosen gluing diffeomorphisms, and is again a holomorphic Lefschetz
fibration. By fiber summing $m$ copies of $f_0$ and $n$ copies of $f_1$ we
obtain a holomorphic Lefschetz fibration with $20m+30n$ singular fibers and
monodromy described by the relation $(\tau_1\tau_2\tau_3\tau_4\tau_5\tau_5
\tau_4\tau_3\tau_2\tau_1)^{2m}(\tau_1\tau_2\tau_3\tau_4\tau_5)^{6n}=1$ in
$\mathrm{Map}_2$. It is worth observing that the fiber sum of
3 copies of $f_0$ is actually the same Lefschetz fibration as the fiber
sum of 2 copies of $f_1$.

It is a recent result of Siebert and Tian \cite{ST} that all genus 2
symplectic Lefschetz fibrations without reducible singular fibers
(i.e.\ without homologically trivial vanishing cycles) and satisfying a
technical assumption (transitivity of monodromy) are holomorphic.
More precisely, there exists a surjective group homomorphism
$\mathrm{Map}_2\to S_6$, mapping each generator $\tau_i$ ($1\le i\le 5$)
to the transposition $(i,i+1)$; the monodromy of a genus 2 Lefschetz
fibration is said to be {\it transitive} if the image of its monodromy
morphism (a subgroup of $\mathrm{Map}_2$) is mapped surjectively onto the
symmetric group $S_6$.

\begin{theorem}[Siebert-Tian]
Every genus 2 Lefschetz fibration without reducible
fibers and with transitive monodromy is holomorphic. It can be
realized as a fiber sum of copies of $f_0$ and $f_1$, and hence its topology
is completely determined by the number of singular fibers.
\end{theorem}

In contrast to this spectacular result, which follows from a positive
solution to the symplectic isotopy problem for certain smooth symplectic
curves in rational ruled surfaces \cite{ST}, if we allow some of the
singular fibers to be reducible, then genus 2 Lefschetz fibrations need not
be holomorphic. Examples of this phenomenon were found by Ozbagci
and Stipsicz \cite{OS}: let $f:X\to S^2$ be the Lefschetz fibration obtained
by blowing up the 4 base points of a pencil of algebraic curves representing
the homology class $[S^2\times pt]+2[pt\times T^2]$ in $S^2\times
T^2$. The Lefschetz fibration $f$ has $8$ singular fibers, two of which are
reducible, and its total space (a blow-up of $S^2\times T^2$) has first
Betti number $b_1=2$. For a suitable choice of the identification
diffeomorphism $\phi$, the symplectic 4-manifold $X\#_\phi X$ obtained by
twisted fiber sum of two copies of the Lefschetz fibration $f$
has first Betti number $b_1=1$, and hence $f\#_\phi f$ is not a holomorphic
Lefschetz fibration.

However, this non-isotopy phenomenon disappears if we perform a
stabilization by fiber sums \cite{Au6}:

\begin{theorem}
Let $f:X\to S^2$ be any genus $2$ Lefschetz fibration. Then the fiber sum of
$f$ with sufficiently many copies of the holomorphic fibration $f_0$
described above is isomorphic to a
holomorphic fibration. Moreover, this fiber sum $f\# n\,f_0$ $(n\gg 0)$
is entirely determined by its total number of singular fibers and by the
number of reducible fibers of each type $($two genus $1$ components, or
genus $0$ and $2$ components$)$.
\end{theorem}

Very little is known about the structure of Lefschetz fibrations with fiber
genus $3$ or more; no analogues of Theorems 3.1 and 3.2 are available.
In fact, by imitating the construction of
Ozbagci and Stipsicz it is easy to construct genus $3$
Lefschetz fibrations without reducible fibers and with $b_1=1$ (which
implies that these fibrations are not holomorphic). More precisely, these
examples can be constructed from the holomorphic genus 3 fibration obtained
by blowing up the 8 base points of a generic pencil of algebraic curves
representing the class $2[S^2\times pt]+2[pt\times T^2]$ in $S^2\times
T^2$. This Lefschetz fibration has $16$ singular fibers, all
irreducible, and its total space has $b_1=2$; suitably twisted fiber sums
of two copies of this fibration yield genus 3 symplectic Lefschetz
fibration structures (without reducible fibers) on symplectic 4-manifolds
with $b_1=1$. Various other examples of non-holomorphic genus 3 Lefschetz
fibrations have been constructed by Amoros et al.\ \cite{ABKP}, Smith
\cite{Sm1}, Fintushel and Stern, among others.

Two natural questions arise at this point. The first one, suggested by the
above examples of non-holomorphic fibrations obtained as twisted fiber sums,
is whether every symplectic Lefschetz fibration can be
decomposed into a fiber sum of holomorphic fibrations. The answer to this
question is known to be negative, as implied by the following result
\cite{Sm2} (see also \cite{ABKP}, \dots):

\begin{theorem}[Amoros et al., Stipsicz, Smith]
If a relatively minimal Lefschetz fibration over $S^2$ with fibers of
genus $\ge 2$ contains a section of square $-1$, then it cannot be decomposed
as any non-trivial fiber sum.
\end{theorem}

In particular, since Lefschetz fibrations obtained by blowing up the base
points of a Lefschetz pencil always admit sections of square $-1$, by applying
Donaldson's construction (Theorem 2.7) to any non-K\"ahler symplectic
4-manifold we obtain indecomposable non-holomorphic Lefschetz fibrations.

The second question we may ask is whether a ``stable isotopy'' result
similar to Theorem 3.2 remains true for Lefschetz fibrations of higher
genus. In the very specific case of Lefschetz fibrations with monodromy
contained in the {\it hyperelliptic} subgroup of the mapping class group,
such a statement can be obtained as a corollary of a recent result of
Kharlamov and Kulikov \cite{KK} about braid monodromy factorizations: 
after repeated (untwisted) fiber sums with copies of a same fixed
holomorphic 
fibration with $8g+4$ singular fibers, any hyperelliptic genus $g$ Lefschetz
fibration eventually becomes holomorphic. Moreover, the fibration obtained
in this manner is completely
determined by its number of singular fibers of each type (irreducible,
reducible with components of given genus), and when the fiber genus is odd
by a certain $\Z_2$-valued invariant. % (spinness if no reducible sings)
(The proof of this result uses the fact that the hyperelliptic mapping
class group is an extension by $\Z_2$ of the braid group of $2g+2$ points
on a sphere, which is itself a quotient of $B_{2g+2}$; this makes it
possible to transform the monodromy of a hyperelliptic Lefschetz fibration
into a factorization in $B_{2g+2}$, with different types of factors
for the various types of singular fibers and extra contributions belonging to
the kernel of the morphism $B_{2g+2}\to B_{2g+2}(S^2)$, and hence reduce the
problem to that studied by Kharlamov and Kulikov). 
However, it is not clear whether the result should be expected to remain
true in the non-hyperelliptic case.

\subsection{Lefschetz fibrations and pseudo-holomorphic curves}
Taubes' work on the Seiberg-Witten invariants of symplectic 4-manifolds
\cite{Ta1,Ta2} has tremendously improved our
understanding of the topology of symplectic 4-manifolds. The most immediate
applications are of two types: to show that certain smooth 4-manifolds admit
no symplectic structure (because their Seiberg-Witten invariants violate
Taubes' structure theorem), or to obtain existence results for
pseudo-holomorphic curves in a given symplectic 4-manifold (using the
non-vanishing of Seiberg-Witten invariants and the identity $SW=Gr_T$).
Results of the second type,
and in particular the fact that any compact symplectic 4-manifold with
$b_2^+\ge 2$ contains an embedded (possibly disconnected) pseudo-holomorphic
curve representing its canonical class \cite{Ta1} (cf.\ Theorem 1.12
$(iv)$), are purely symplectic in nature, and it may be possible to reprove
them by methods that do not involve Seiberg-Witten theory at all.
Spectacular progress has been made in this direction by Donaldson and
Smith, who have obtained the following result \cite{DS}:

\begin{theorem}[Taubes, Donaldson-Smith]
Let $(X,\omega)$ be a compact symplectic 4-manifold with
$b_2^+(X)>1+b_1(X)$, and assume that the cohomology class $[\omega]\in
H^2(X,\R)$ is rational. Then there exists a smooth (not necessarily
connected) embedded symplectic
surface in $X$ which represents the homology class Poincar\'e dual to
$c_1(K_X)$.
\end{theorem}

Following the work of Donaldson and Smith \cite{DS,Sm3}, we now give an
outline of the proof of this result. The main idea is to represent a blow-up
$\hat{X}$ of
the manifold $X$ as a symplectic Lefschetz fibration, and to look for
symplectic submanifolds embedded in ``standard'' position in $\hat{X}$
with respect to the Lefschetz fibration $f:\hat{X}\to S^2$. Namely,
$\Sigma\subset \hat{X}$ is said to be a {\it standard surface} if it
avoids the critical points of $f$ and if it is positively
transverse to the fibers of $f$ everywhere except at isolated
non-degenerate tangency points (i.e., the restriction of $f$ to $\Sigma$
defines a branched covering with simple branch points).

Let $r>0$ be the intersection number of $[\Sigma]$ with the fiber of $f$
(in the case of the canonical class,
$r=2g-2$, where $g$ is the fiber genus), and consider a fibration 
$F:X_r(f)\to S^2$ with generic fiber the $r$-th symmetric product of the
fiber of $f$. More precisely, choosing an almost-complex structure in order
to make each fiber of $f$ a (possibly nodal) Riemann surface, $X_r(f)$ is
defined as the {\it relative Hilbert scheme} of degree $r$ divisors in the
fibers of $f$; each fiber of $F$ corresponding to a smooth
fiber of $f$ is naturally identified with its $r$-fold symmetric product,
but in the case of nodal fibers the Hilbert scheme is a partial
desingularization of the symmetric product.

A standard surface $\Sigma$ intersects every fiber of $f$ in $r$
points (counting with multiplicities), which defines a point in the $r$-th
symmetric product (or Hilbert scheme) of the fiber. Hence, to a standard
surface $\Sigma\subset\hat{X}$ we can associate a section $s_\Sigma$ of the
fibration $F:X_r(f)\to
S^2$.  The points where $\Sigma$ becomes non-degenerately tangent to
the fibers of $f$ correspond to positive transverse intersections between
the section $s_\Sigma$ and the diagonal divisor $\Delta\subset X_r(f)$
consisting of all $r$-tuples in which two or more points coincide.
Conversely, any section of $X_r(f)$ that intersects $\Delta$ transversely
and positively determines a standard surface in $\hat{X}$.

There is a well-defined map from the space of homotopy classes of sections
of $X_r(f)$ to the homology group $H_2(\hat{X},\Z)$, which to a section of
$X_r(f)$ associates the homology class represented by the corresponding
surface in $\hat{X}$. This map is injective, so to a given homology class
$\alpha\in H_2(X',\Z)$ we can associate at most a single homotopy class
$\bar{\alpha}$ of sections of $X_r(f)$ ($r=\alpha\cdot [f^{-1}(pt)]$).

We can equip the manifold $X_r(f)$ with symplectic and almost-complex
structures whose restrictions to each fiber of $F$ coincide with
the standard K\"ahler
and complex structures on the symmetric product of a Riemann surface.
Then it is possible to define an invariant $Gr_{DS}(X,f;\alpha)\in\Z$
counting pseudo-holomorphic sections of $X_r(f)$ in the homotopy class
corresponding to a given homology class $\alpha\in H_2(X',\Z)$. More
precisely, if the expected dimension of the moduli space of
pseudo-holomorphic sections is positive then one obtains an integer-valued
invariant by adding incidence conditions requiring the standard surface
$\Sigma$ to pass through certain points in $X'$, or equivalently requiring
the section $s_\Sigma$ to pass through certain divisors in the corresponding
fibers of $F$. The fact that we are interested in a moduli space of
pseudo-holomorphic sections makes it possible to control bubbling,
which can only occur inside the fibers of $F$.

In the case of the canonical class $\alpha=c_1(K_{\hat{X}})$, there exists a
specific almost-complex structure on $X_{2g-2}(f)$ for which the moduli space
of pseudo-holomorphic sections in the relevant homotopy class $\bar\alpha$ is
a projective space of complex dimension $\frac{1}{2}(b_2^+(X)-1-b_1(X))-1$
(whereas the generically expected dimension is $0$). By
computing the Euler class of the obstruction bundle over this moduli space,
one obtains that $Gr_{DS}(X,f;c_1(K_{\hat{X}}))=\pm 1$ \cite{DS}.

The non-vanishing of this invariant gives us an existence result for
pseudo-holomorphic sections of $X_r(f)$. However, in order to obtain a
standard symplectic surface from such a section $s$, we need to ensure the
positivity of the intersections between $s$ and the diagonal divisor
$\Delta$. This means that we need to consider on $X_r(f)$ a different
almost-complex structure, for which the strata composing $\Delta$ are
pseudo-holomorphic submanifolds of $X_r(f)$ (as well as the divisors
associated with the exceptional sections of $f$ arising from the blow-ups
of the base points of the chosen pencil on $X$). For a generic choice of
almost-complex structure compatible with the diagonal strata, the moduli
spaces of pseudo-holomorphic sections of $F$ have the expected dimension,
and the non-vanishing of the $Gr_{DS}$ invariant leads to the existence of
a section transverse to all diagonal strata in which it is not contained.

From such a section of $F$ we can obtain a (possibly disconnected)
standard symplectic surface $\hat\Sigma\subset\hat{X}$, representing
the homology class $c_1(K_{\hat{X}})=\pi^* c_1(K_X)+
\sum [E_i]$, where $E_i$ are the exceptional spheres of the blow-up
$\pi:\hat{X}\to X$. Moreover, because this surface has intersection number
$-1$ with each sphere $E_i$, the local positivity of intersection implies
that $E_i$ is actually contained in $\hat\Sigma$, and therefore we have
$\hat\Sigma=\Sigma\cup\bigcup E_i$, for some symplectic surface
$\Sigma\subset\hat{X}-\bigcup E_i$; the image of $\Sigma$ in $X$ is a smooth
symplectic submanifold representing the homology class Poincar\'e dual to
the canonical bundle \cite{DS}.

While it has not yet been shown that the invariant $Gr_{DS}$ is actually
independent of the choice of a high degree Lefschetz pencil structure on
$(X^4,\omega)$ and coincides with the invariant $Gr_T$ defined by Taubes,
it is worth mentioning a Serre duality-type result obtained by Smith for
the $Gr_{DS}$ invariant: under the same assumption as above on the symplectic
4-manifold $X$ ($b_2^+(X)>1+b_1(X)$), for a Lefschetz pencil of sufficiently
high degree we have $Gr_{DS}(X,f;\alpha)=\pm Gr_{DS}(X,f;K_X-\alpha)$
\cite{Sm3}. This is to be compared to the duality formula for Gromov-Taubes
invariants, $Gr_T(\alpha)=\pm Gr_T(K_X-\alpha)$,
which follows immediately from a duality among Seiberg-Witten invariants
($SW(-L)=\pm SW(L)$) and Taubes' result (Theorem~1.12).

\subsection{Fukaya-Seidel categories for Lefschetz pencils}
One of the most exciting applications of Lefschetz pencils, closely
related to the homological mirror symmetry conjecture, is Seidel's
construction of ``directed Fukaya categories'' associated to a Lefschetz
pencil. Besides the many technical difficulties arising in their definition,
Fukaya categories of symplectic manifolds are intrinsically very hard to
compute, because relatively little is known about embedded Lagrangian
submanifolds in symplectic manifolds of dimension $\ge 4$, especially in
comparison to the much better understood theory of coherent sheaves over
complex varieties, which play the role of their mirror counterparts.

Consider an arc $\gamma$ joining a regular value $p_0$ to a critical value
$p_1$ in the base of a symplectic Lefschetz fibration $f:X\to S^2$.
Using the horizontal distribution given by the symplectic orthogonal to
the fibers, we can transport the vanishing cycle at $p_1$ along the arc
$\gamma$ to obtain a Lagrangian disc $D_\gamma\subset X$ fibered above $\gamma$,
whose boundary is an embedded Lagrangian sphere $S_\gamma$ in the fiber
$\Sigma_0=f^{-1}(p_0)$. The Lagrangian disc $D_\gamma$ is sometimes called
the {\it Lefschetz thimble} over $\gamma$, and its boundary $S_\gamma$ is the
vanishing cycle already considered in \S 2. If we now consider an
arc $\gamma$ joining two critical values $p_1,p_2$ of $f$ and passing
through $p_0$, then the above construction applied to each half of $\gamma$
yields two Lefschetz thimbles $D_1$ and $D_2$, whose boundaries are
Lagrangian spheres $S_1,S_2\subset \Sigma_0$. If $S_1$ and $S_2$ coincide
exactly, then $D_1\cup D_2$ is an embedded Lagrangian sphere in $X$,
fibering above the arc $\gamma$ (see the picture below); more generally,
if $S_1$ and $S_2$ are Hamiltonian isotopic to each other, then perturbing
slightly the symplectic structure we can reduce to the previous case and
obtain again a Lagrangian sphere in $X$. The arc
$\gamma$ is called a {\it matching path} in the Lefschetz fibration $f$.
\medskip

\setlength{\unitlength}{0.5mm}
\begin{center}
\begin{picture}(50,55)(-10,-15)
\put(-10,0){\line(1,0){50}}
\put(-10,40){\line(1,0){50}}
\qbezier[30](30,40)(35,40)(35,35)
\qbezier[30](35,35)(35,30)(34,27.5)
\qbezier[30](34,27.5)(33,25)(33,20)
\qbezier[30](30,0)(35,0)(35,5)
\qbezier[30](35,5)(35,10)(34,12.5)
\qbezier[30](34,12.5)(33,15)(33,20)
\qbezier[30](30,40)(25,40)(25,36)
\qbezier[30](25,36)(25,32)(26,28)
\qbezier[20](26,28)(27,26)(27,25)
\qbezier[30](30,0)(25,0)(25,5)
\qbezier[30](25,5)(25,10)(26,12.5)
\qbezier[80](26,12.5)(27,15)(27,25)
\qbezier[30](32,32)(32,36)(30,36)
\qbezier[30](30,36)(28,36)(26,32)
\qbezier[20](32,32)(32,30)(30,30)
\qbezier[20](30,30)(28,30)(26,32)
\qbezier[30](29,5)(33,8)(29,11)
\qbezier[20](30,6)(28,8)(30,10)
\qbezier[30](0,40)(5,40)(5,35)
\qbezier[30](5,35)(5,30)(4,27.5)
\qbezier[30](4,27.5)(3,25)(3,20)
\qbezier[30](0,0)(5,0)(5,5)
\qbezier[30](5,5)(5,10)(4,12.5)
\qbezier[30](4,12.5)(3,15)(3,20)
\qbezier[30](0,40)(-5,40)(-5,36)
\qbezier[30](-5,36)(-5,32)(-4,28)
\qbezier[20](-4,28)(-3,26)(-3,25)
\qbezier[30](0,0)(-5,0)(-5,5)
\qbezier[30](-5,5)(-5,10)(-4,12.5)
\qbezier[80](-4,12.5)(-3,15)(-3,25)
\qbezier[30](2,32)(2,36)(0,36)
\qbezier[30](0,36)(-2,36)(-4,32)
\qbezier[20](2,32)(2,30)(0,30)
\qbezier[20](0,30)(-2,30)(-4,32)
\qbezier[30](-1,5)(3,8)(-1,11)
\qbezier[20](0,6)(-2,8)(0,10)
\put(-4,32){\circle*{1.5}}
\put(26,32){\circle*{1.5}}
\qbezier[80](-4,32)(-4,36)(11,36)
\qbezier[80](-4,32)(-4,28)(11,28)
\qbezier[80](26,32)(26,28)(11,28)
\qbezier[80](26,32)(26,36)(11,36)
\qbezier[30](11,28)(9,32)(11,36)
\qbezier[10](11,28)(13,32)(11,36)
\put(9.5,20){$S^n$}
\put(-12,-15){\line(1,0){44}}
\put(-12,-15){\line(1,1){8}}
\put(32,-15){\line(1,1){8}}
\put(0,-10){\circle*{1.5}}
\put(30,-10){\circle*{1.5}}
\put(0,-10){\line(1,0){30}}
\put(15,-7){$\gamma$}
\end{picture}
\end{center}

\noindent
Matching paths are an important source of Lagrangian spheres, and more
generally (extending suitably the notion of matching path) of embedded
Lagrangian submanifolds. Conversely, a folklore theorem asserts that
any given embedded Lagrangian sphere (or more generally, compact Lagrangian
submanifold) in a compact symplectic manifold is isotopic to one that fibers
above a matching path in a Donaldson-type symplectic Lefschetz pencil of
sufficiently high degree.

The intersection theory of Lagrangian spheres that fiber above matching
paths is much nicer than that of arbitrary Lagrangian spheres, because
if two Lagrangian spheres $S,S'\subset X$ fiber above matching paths
$\gamma,\gamma'$, then all intersections of $S$ with $S'$ lie in the
fibers above the intersection points of $\gamma$ with $\gamma'$; hence,
the Floer homology of $S$ and $S'$ can be computed by studying intersection
theory for Lagrangian spheres in the fibers of $f$, rather than in its
total space. 

These considerations have led Seidel to the following construction of a
Fukaya-type $A_\infty$-category associated to a symplectic Lefschetz
pencil $f$ on a compact symplectic manifold $(X,\omega)$ \cite{Se1}.
Let $f$ be a symplectic Lefschetz pencil determined by two sections
$s_0,s_1$ of a sufficiently positive line bundle $L^{\otimes k}$ as in Theorem 2.7.
Assume that $\Sigma_\infty=s_1^{-1}(0)$ is a smooth fiber of the pencil,
and consider the symplectic manifold with boundary $X^0$ obtained from $X$
by removing a suitable neighborhood of $\Sigma_\infty$. The map $f$ induces
a symplectic Lefschetz fibration structure $f^0:X^0\to D^2$ over a disc,
whose fibers are symplectic submanifolds with boundary obtained from the
fibers of $f$ by removing a neighborhood of their intersection points with
the symplectic hypersurface $\Sigma_\infty$ (the base points of the pencil).
Choose a reference point $p_0\in \partial D^2$, and consider the fiber
$\Sigma_0=(f^0)^{-1}(p_0)\subset X^0$.

Let $\gamma_1,\dots,\gamma_r$ be a collection of arcs in $D^2$ joining the
reference point $p_0$ to the various critical values of $f^0$, intersecting
each other only at $p_0$, and ordered in the clockwise direction around
$p_0$. As discussed above, each arc $\gamma_i$ gives rise to a Lefschetz
thimble $D_i\subset X^0$, whose boundary is a Lagrangian sphere $L_i\subset
\Sigma_0$. To avoid having to discuss the orientation of moduli spaces, we
give the following definition using $\Z_2$ as the coefficient ring
\cite{Se1}:

\begin{definition}[Seidel]
The directed Fukaya category $FS(f;\{\gamma_i\})$ is the following
$A_\infty$-category: the objects of $FS(f;\{\gamma_i\})$ are the Lagrangian
vanishing cycles $L_1,\dots,L_r$; the morphisms between
the objects are given by
$$\mathrm{Hom}(L_i,L_j)=\begin{cases}
CF^*(L_i,L_j;\Z_2)=\Z_2^{|L_i\cap L_j|} & \mathrm{if}\ i<j\\
\Z_2\,\,id & \mathrm{if}\ i=j\\
0 & \mathrm{if}\ i>j;
\end{cases}$$
and the differential $\mu^1$, composition $\mu^2$ and higher order
compositions $\mu^n$ are given by Lagrangian Floer homology inside
$\Sigma_0$: more precisely,
$$\mu^n:\mathrm{Hom}(L_{i_0},L_{i_1})\otimes \dots\otimes
\mathrm{Hom}(L_{i_{n-1}},L_{i_n}) \to \mathrm{Hom}(L_{i_0},L_{i_n})[2-n]$$
is trivial when the inequality $i_0<i_1<\dots<i_n$ fails to hold (i.e.\ it
is always zero in this case, except for $\mu^2$ where composition with an
identity morphism is given by the obvious formula).
When $i_0<\dots<i_n$, $\mu^n$ is defined
by counting pseudo-holomorphic maps from the disc to $\Sigma_0$, mapping
$n$ cyclically ordered marked points on the boundary to the given
intersection points between vanishing cycles, and the portions of boundary
between them to $L_{i_0},\dots,L_{i_n}$.
\end{definition}

One of the most attractive features of this definition is that it only
involves Floer homology for Lagrangians inside the hypersurface $\Sigma_0$;
in particular, when $X$ is a symplectic 4-manifold, the definition becomes
essentially combinatorial, since in the case of a Riemann surface
the pseudo-holomorphic discs appearing in the definition of Floer homology
and product structures are essentially topological objects.

From a technical point of view, a key property that greatly facilitates the
definition of Floer homology for the vanishing cycles $L_i$ is {\it
exactness}. Namely, the symplectic structure on the manifold
$X^0$ is {\it exact}, i.e.\ it can
be expressed as $\omega=d\theta$ for some 1-form $\theta$ (up to a scaling
factor, $\theta$ is the 1-form describing the connection on $L^{\otimes k}$
in the trivialization of $L^{\otimes k}$ over $X-\Sigma_\infty$ induced by
the section $s_1/|s_1|$). With this understood, the submanifolds
$L_i$ are all {\it exact Lagrangian}, i.e.\ the restriction $\theta_{|L_i}$
is not only closed ($d\theta_{|L_i}=\omega_{|L_i}=0$) but also exact,
$\theta_{|L_i}=d\phi_i$. Exactness has two particularly nice consequences.
First, $\Sigma^0$ contains no closed pseudo-holomorphic curves (because
the cohomology class of $\omega=d\theta$ vanishes). Secondly, there are no
non-trivial pseudo-holomorphic discs in $\Sigma_0$ with boundary contained
in one of the Lagrangian submanifolds $L_i$.
Indeed, for any such disc $D$, we have $\mathrm{Area}(D)=\int_D\omega=
\int_{\partial D}\theta=\int_{\partial D} d\phi_i=0$. Therefore, bubbling
never occurs (neither in the interior nor on the boundary of the domain) in
the moduli spaces used to define the Floer homology groups
$HF(L_i,L_j)$. Moreover, the exactness of $L_i$ provides a priori estimates
on the area of all pseudo-holomorphic discs contributing to the definition
of the products $\mu^n$ ($n\ge 1$); this implies the
finiteness of the number of discs to be considered and solves elegantly the
convergence problems that normally make it necessary to define Floer
homology over Novikov rings.

It is natural question to ask oneself how much the category
$FS(f,\{\gamma_i\})$ depends on the given data, and in particular on the
chosen ordered collection of arcs $\{\gamma_i\}$. An answer is provided by a
result of Seidel showing that, if the ordered collection $\{\gamma_i\}$ is
replaced by another one $\{\gamma'_i\}$, then the categories
$FS(f,\{\gamma_i\})$ and $FS(f,\{\gamma'_i\})$ differ by a sequence of
{\it mutations} (operations that modify the ordering of the objects of
the category while twisting some of them along others) \cite{Se1}.
Hence, the category naturally associated to the Lefschetz
pencil $f$ is not the finite directed category defined above, but rather a
{\it derived} category, obtained from $FS(f,\{\gamma_i\})$ by considering
formal direct sums and twisted complexes of Lagrangian vanishing cycles
(with additional features such as idempotent splittings, formal inverses
of quasi-isomorphisms, ...). It is a classical result that, if two
categories differ by mutations, then their derived categories are
equivalent; hence the derived category $D(FS(f))$ only depends on the
Lefschetz pencil $f$ rather than on the choice of an ordered system of arcs
\cite{Se1}.

As an example, let us consider the case of a pencil of conics on $\CP^2$,
an example for which monodromy was described explicitly at the beginning of
\S 3.1. The fiber $\Sigma_0$ is a sphere with four punctures, and for a
suitable system of arcs the three vanishing cycles
$L_1,L_2,L_3\subset\Sigma_0$ are as represented at the beginning of \S 3.1:
\medskip

\begin{center}
\epsfig{file=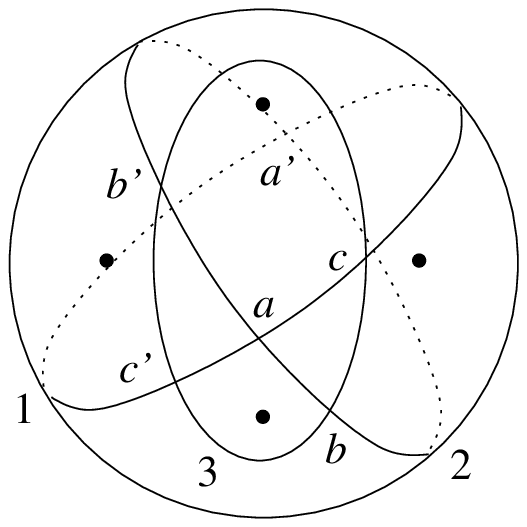,height=1.2in}
\end{center}

Any two of these three vanishing cycles intersect transversely in two
points, so $\mathrm{Hom}(L_1,L_2)=\Z_2\,a\oplus \Z_2\,a'$,
$\mathrm{Hom}(L_2,L_3)=\Z_2\,b\oplus \Z_2\,b'$, and
$\mathrm{Hom}(L_1,L_3)=\Z_2\,c\oplus \Z_2\,c'$ are all two-dimensional.
There are no embedded 2-sided polygons in the punctured sphere $\Sigma_0$
with boundary in $L_i\cup L_j$ for any pair $(i,j)$, since each of the four
regions delimited by $L_i$ and $L_j$ contains one of the punctures, so
$\mu^1\equiv 0$. However, there are four triangles with boundary in
$L_1\cup L_2\cup L_3$ (with vertices $abc$, $ab'c'$, $a'b'c$, $a'bc'$
respectively), and in each case the cyclic ordering of the boundary is
compatible with the ordering of the vanishing cycles. Therefore, the
composition of morphisms is given by the formulas
$\mu^2(a,b)=\mu^2(a',b')=c$, $\mu^2(a,b')=\mu^2(a',b)=c'$.
Finally, the higher compositions $\mu^n$, $n\ge 3$ are all trivial in this
category, because the ordering condition $i_0<\dots<i_n$ never holds~\cite{Se2}.

We finish this section by some vague considerations about the relationship
between the derived Fukaya-Seidel category $D(FS(f))$ of the Lefschetz
pencil $f$ and the derived Fukaya category $D\mathcal{F}(X)$ of the
symplectic manifold $X$, following ideas of Seidel. 
At first glance, the category $FS(f;\{\gamma_i\})$
appears to be more closely related to the Fukaya category of the fiber
$\Sigma_0$ than to that of the total space of the Lefschetz pencil. However,
the objects $L_i$ actually correspond not only to Lagrangian spheres in
$\Sigma_0$ (the vanishing cycles), but also to Lagrangian discs in $X^0$
(the Lefschetz thimbles $D_i$); and the Floer intersection theory in
$\Sigma_0$ giving rise to $\mathrm{Hom}(L_i,L_j)$ and to the product
structures can also be thought of in terms of intersection theory for
the Lagrangian discs $D_i$ in $X^0$. When passing to the derived category
$D(FS(f))$, we hugely increase the number of objects, by considering not
only the thimbles $D_i$ but also arbitrary complexes obtained from them;
this means that the objects of $D(FS(f))$ include arbitrary
(not necessarily closed) Lagrangian submanifolds in $X^0$, with boundary
in $\Sigma_0$. Since Fukaya categories are only concerned with closed
Lagrangian submanifolds, it is necessary to consider a subcategory of
$D(FS(f))$ whose objects correspond only to the closed Lagrangian
submanifolds in $X^0$ (i.e., combinations of $D_i$ for which the
boundaries cancel); it is expected that this can be done in purely
categorical terms by considering those objects of $D(FS(f))$ on which
the Serre functor acts simply by a shift. The resulting subcategory
should be closely related to the derived Fukaya category of the open
manifold $X^0$. This leaves us with the problem of relating
$\mathcal{F}(X^0)$ with $\mathcal{F}(X)$. These two categories have
the same objects and morphisms (Lagrangians in $X$ can be made disjoint
from $\Sigma_\infty$), but the differentials and product structures differ.
More precisely, the definition of $\mu^n$ in $\mathcal{F}(X^0)$ only
involves counting pseudo-holomorphic discs contained in $X^0$, i.e.\ disjoint
from the hypersurface $\Sigma_\infty$. In order to account for the missing
contributions, one should introduce a formal parameter $q$ and count the
pseudo-holomorphic discs with boundary in $\bigcup L_i$ that intersect
$\Sigma_\infty$ in $m$ points (with multiplicities) with a coefficient
$q^m$. The introduction of this parameter $q$ leads to a {\it deformation of
$A_\infty$-structures}, i.e.\ an
$A_\infty$-category in which the differentials and products $\mu^n$ are
defined over a ring of
formal power series in the variable $q$; the limit $q=0$ corresponds to
the (derived) Fukaya category $D\mathcal{F}(X^0)$, while non-zero values of
$q$ are expected to yield $D\mathcal{F}(X)$.

The above considerations provide a strategy that should make it possible (at
least in some examples) to calculate Fukaya categories by induction on
dimension; an important consequence is that it becomes possible to verify
the homological mirror symmetry conjecture (or parts thereof) on a wider
class of examples (e.g.\ some K3 surfaces, cf.\ recent work of Seidel).
%}

\section{Symplectic branched covers of $\CP^2$}

%\comment{

\subsection{Symplectic branched covers}

\begin{definition}
A smooth map $f:X^4\to (Y^4,\omega_Y)$ from a compact oriented smooth
4-manifold to a compact symplectic 4-manifold is a {\em symplectic branched
covering} if, given any point $p\in X$, there exist 
neighborhoods $U\ni p$ and $V\ni f(p)$ and orientation-preserving local
diffeomorphisms $\phi:U\to\C^2$ and $\psi:V\to\C^2$, such that
$\psi_*\omega_Y(v,iv)>0$ $\forall v\neq 0$ $($i.e.\ the standard complex
structure is $\psi_*\omega_Y$-tame$)$, and such that
$\psi\circ f\circ \phi^{-1}$ is one of the following model maps:

$(i)$ $(u,v)\mapsto (u,v)$ (local diffeomorphism),

$(ii)$ $(u,v)\mapsto (u^2,v)$ (simple branching),

$(iii)$ $(u,v)\mapsto (u^3-uv,v)$ (cusp).
\end{definition}

The three local models appearing in this definition are exactly those
describing a generic holomorphic map between complex surfaces, except that
the local coordinate systems we consider are not holomorphic.

By computing the Jacobian of $f$ in the given local coordinates, we can see
that the {\it ramification curve} $R\subset X$ is a smooth submanifold (it
is given by $\{u=0\}$ in the second local model and $\{v=3u^2\}$ in the
third one). However, the image $D=f(R)\subset X$ (the {\it branch curve},
or {\it discriminant curve}) may be singular. More precisely, in the simple
branching model $D$ is given by $\{z_1=0\}$, while in the cusp model we have
$f(u,3u^2)=(-2u^3,3u^2)$, and hence $D$ is locally identified with
the singular curve $\{27z_1^2=4z_2^3\}\subset\C^2$. This means that, at the
cusp points, $D$ fails to be immersed. Besides the cusps, the branch curve
$D$ also generically presents {\it transverse double points} (or {\it nodes}),
which do not appear in the local models because they correspond to simple
branching in two distinct points $p_1,p_2$ of the same fibre of $f$. 
There is no constraint on the orientation of the local intersection between
the the two branches of $D$ at a node (positive or negative,
i.e.\ complex or anti-complex), because the local models near $p_1$ and
$p_2$ hold in different coordinate systems on $Y$.

Generically, the only singularities of the branch curve
$D\subset Y$ are transverse double points (``nodes'') of either orientation
and complex cusps. Moreover, because the local
models identify $D$ with a complex curve, the tameness condition on the
coordinate systems implies that $D$ is a (singular) symplectic submanifold
of~$Y$.

The following result states that a symplectic branched cover of a symplectic
4-manifold carries a natural symplectic structure \cite{Au2}:

\begin{proposition}
If $f:X^4\to (Y^4,\omega_Y)$ is a symplectic branched cover, then $X$
carries a symplectic form $\omega_X$ such that $[\omega_X]=f^*[\omega_Y]$,
canonically determined up to symplectomorphism.
\end{proposition}

\proof
The 2-form $f^*\omega_Y$ is closed, but it is only non-degenerate outside of
$R$. At any point $p$ of $R$, the 2-plane $K_p=\mathrm{Ker}\,df_p\subset
T_pX$ carries a natural orientation induced by the complex orientation in
the local coordinates of Definition 4.1. Using the local models, we can
construct an {\it exact} 2-form $\alpha$ such that, at any point $p\in R$,
the restriction of $\alpha$ to $K_p$ is non-degenerate and positive. 

More precisely, given $p\in R$ we consider a small ball centered at $p$ and 
local coordinates $(u,v)$ such that $f$ is given by one of the models of
Definition 4.1, and we set $\alpha_p=d(\chi_1(|u|)\chi_2(|v|)\, x\,dy)$,
where $x=\mathrm{Re}(u)$, $y=\mathrm{Im}(u)$, and $\chi_1$ and $\chi_2$ are
suitably chosen smooth cut-off functions. We then define $\alpha$ to be
the sum of these $\alpha_p$ when $p$ ranges over a
finite subset of $R$ for which the supports of the $\alpha_p$ cover the
entire ramification curve $R$. 
Since $f^*\omega_Y\wedge \alpha$ is positive at every point of $R$, it
is easy to check that the 2-form $\omega_X=f^*\omega_Y + \epsilon\,\alpha$
is symplectic for a small enough value of the constant $\epsilon>0$.

The fact that $\omega_X$ is canonical up to symplectomorphism follows
immediately from Moser's stability theorem and from the observation that
the space of exact perturbations
$\alpha$ such that $\alpha_{|K_p}>0$ $\forall p\in R$ is a convex subset
of $\Omega^2(X)$ and hence connected.
\endproof

\subsection{Symplectic manifolds and maps to $\CP^2$}

Approximately holomorphic techniques make it possible to show that every
compact symplectic 4-manifold can be realized as a branched cover of
$\CP^2$. The general setup is
similar to Donaldson's construction of symplectic Lefschetz pencils: we
consider a compact symplectic manifold $(X,\omega)$, and perturbing the
symplectic structure if necessary we may assume that
$\frac{1}{2\pi}[\omega]\in H^2(X,\Z)$. Introducing an almost-complex 
structure $J$ and a line bundle $L$ with $c_1(L)=\frac{1}{2\pi}[\omega]$,
we consider triples of approximately holomorphic sections
$(s_k^0,s_k^1,s_k^2)$ of $L^{\otimes k}$: for $k\gg 0$, it is again possible
to achieve a generic behavior for the projective map
$f_k=(s_k^0:s_k^1:s_k^2):X\to\CP^2$ associated with the
linear system. If the manifold $X$ is four-dimensional, then the linear
system generically has no base points, and for a suitable choice of sections
the map $f_k$ is a branched covering \cite{Au2}.

\begin{theorem}
For large enough $k$, three suitably chosen approximately holomorphic
sections of $L^{\otimes k}$ over $(X^4,\omega)$ determine a symplectic
branched covering $f_k:X^4\to\CP^2$, described in approximately holomorphic
local coordinates by the local models of Definition 4.1. Moreover, for $k\gg
0$ these branched covering structures are uniquely determined up to isotopy.
\end{theorem}

Because the local models hold in approximately holomorphic (and hence
$\omega$-tame) coordinates, the ramification curve $R_k$ of $f_k$ is a
symplectic submanifold in $X$ (connected, since the Lefschetz hyperplane
theorem applies). Moreover, if we normalize the Fubini-Study symplectic form
on $\CP^2$ in such a way that $\frac{1}{2\pi}[\omega_{FS}]$ is the generator
of $H^2(\CP^2,\Z)$, then we have $[f_k^*\omega_{FS}]=2\pi c_1(L^{\otimes k}=
k[\omega]$, and it is fairly easy to check that the symplectic form on $X$
obtained by applying Proposition 4.2 to the branched covering $f_k$
coincides up to symplectomorphism with $k\omega$ \cite{Au2}. In fact, the
exact 2-form
$\alpha=k\omega-f_k^*\omega_{FS}$ is positive over $\mathrm{Ker}\,df_k$ at
every point of $R_k$, and $f_k^*\omega_{FS}+t\alpha$ is a symplectic form
for all $t\in (0,1]$.

The uniqueness statement in Theorem 4.3, which should be interpreted exactly
in the same way as that obtained by Donaldson for Lefschetz pencils, implies
that for $k\gg 0$ it is possible to define invariants of the symplectic
manifold $(X,\omega)$ in terms of the monodromy of the branched covering
$f_k$ and the topology of its branch curve $D_k\subset \CP^2$. However,
the branch curve $D_k$ is only determined up to creation or cancellation
of (admissible) pairs of nodes of opposite orientations.

A similar construction can be attempted when $\dim X>4$; in this case, the
set of base points $Z_k=\{s_k^0=s_k^1=s_k^2=0\}$ is no longer empty. The set
of base points is generically a smooth codimension 6 symplectic submanifold.
With this understood, Theorem 4.3 admits the following higher-dimensional
analogue \cite{Au3}:

\begin{theorem}
For large enough $k$, three suitably chosen approximately holomorphic
sections of $L^{\otimes k}$ over $(X^{2n},\omega)$ determine a map
$f_k:X-Z_k\to\CP^2$ with generic local models, canonically determined up to
isotopy.
\end{theorem}

The model maps describing the local behavior of $f_k$ in approximately
holomorphic local coordinates are now the following:\smallskip

\begin{tabular}{cl}
$(0)$ & $(z_1,\dots,z_n)\mapsto (z_1:z_2:z_3)$ near a base point,\\
$(i)$ & $(z_1,\dots,z_n)\mapsto (z_1,z_2)$,\\
$(ii)$ & $(z_1,\dots,z_n)\mapsto (z_1^2+\dots+z_{n-1}^2,z_n)$,\\
$(iii)$ & $(z_1,\dots,z_n)\mapsto (z_1^3-z_1z_n+z_2^2+\dots+z_{n-1}^2,z_n)$.
\end{tabular}
\smallskip

\noindent
The set of critical points $R_k\subset X$ is again a (connected) smooth
symplectic curve, and its image $D_k=f_k(R_k)\subset\CP^2$ is again a
singular symplectic curve whose only singularities generically are
transverse double points of either orientation and complex cusps.
The fibers of $f_k$ are codimension 4 symplectic submanifolds, intersecting
along $Z_k$; the fiber above a point of $\CP^2-D_k$ is smooth, while
the fiber above a smooth point of $D_k$ presents an ordinary double point,
the fiber above a node presents two ordinary double points,
and the fiber above a cusp presents an $\mathrm{A}_2$ singularity.

As in the four-dimensional case, the uniqueness statement implies that, up
to possible creations or cancellations of pairs of double points with
opposite orientations in the curve $D_k$, the topology of the fibration
$f_k$ can be used to define invariants of the manifold
$(X^{2n},\omega)$.

The proof of Theorems 4.3 and 4.4 relies on a careful examination of the
various possible local behaviors for the map $f_k$ and on transversality
arguments showing the existence of sections of $L^{\otimes k}$ with generic
behavior. Hence, the argument relies on the enumeration of the various
special cases, generic or not, that may occur; each one corresponds to the
vanishing of a certain quantity that can be expressed in terms of the
sections $s_k^0,s_k^1,s_k^2$ and their derivatives.

In order to simplify these arguments, and to make it possible to
extend these results to linear systems generated by more than three sections
or even more general situations, it is helpful to develop an approximately
holomorphic version of singularity theory. The core ingredient of this
approach is a uniform transversality result for jets of approximately
holomorphic sections \cite{Au4}.

Given approximately holomorphic sections $s_k$ of very positive bundles
$E_k$ (e.g.\ $E_k=\C^m\otimes L^{\otimes k}$) over the symplectic manifold
$X$, one can consider the $r$-jets $j^r s_k=(s_k,\partial s_k,(\partial
\partial s_k)_\mathrm{sym}, \dots, (\partial^r s_k)_\mathrm{sym})$, which
are sections of the {\it jet bundles} $\mathcal{J}^r E_k=\bigoplus_{j=0}^r
(T^*X^{(1,0)})_{\mathrm{sym}}^{\otimes j}\otimes E_k$. Jet bundles can
naturally be stratified by approximately holomorphic submanifolds
corresponding to the various possible local behaviors at order $r$ for the
sections $s_k$. The generically expected behavior corresponds to the case
where the jet $j^r s_k$ is transerse to the submanifolds in the
stratification. The result is the following \cite{Au4}:

\begin{theorem}
Given stratifications $\mathcal{S}_k$ of the jet bundles $\mathcal{J}^r E_k$
by a finite number of approximately holomorphic submanifolds
(Whitney-regular, uniformly transverse to fibers, and with curvature bounded
independently of $k$), for large enough $k$ the vector bundles $E_k$ admit
approximately holomorphic sections $s_k$ whose $r$-jets are uniformly
transverse to the stratifications $\mathcal{S}_k$. Moreover these sections
may be chosen arbitrarily close to given sections.
\end{theorem}

A one-parameter version of this result also holds, which makes it
possible to obtain results of asymptotic uniqueness up to isotopy for
generic sections~\cite{Au4}.

Applied to suitably chosen stratifications, Theorem 4.5 provides the main
ingredient for the construction of $m$-tuples of approximately
holomorphic sections of $L^{\otimes k}$ (and hence projective maps $f_k$ with
values in $\CP^{m-1}$) with generic behavior. Once uniform transversality
of jets has been obtained, the only remaining task is to achieve some
control over the antiholomorphic derivative $\dbar f_k$ near the critical
points of $f_k$ (typically its vanishing in some directions), in order to
ensure that $\dbar f_k \ll \partial f_k$ everywhere; for low values of $m$
such as those considered above, this task is comparatively easy.

\subsection{Monodromy invariants for branched covers of $\CP^2$}

The topological data characterizing a symplectic branched covering
$f:X^4\to\CP^2$ are on one hand the topology of the branch curve
$D\subset\CP^2$ (up to isotopy and cancellation of pairs of nodes),
and on the other hand a monodromy morphism $\theta:\pi_1(\CP^2-D)\to S_N$
describing the manner in which the $N=\deg f$ sheets of the covering are
arranged above $\CP^2-D$.

Some simple properties of the monodromy morphism $\theta$ can be readily
seen by considering the local models of Definition 4.1. For example, the
image of a small loop $\gamma$ bounding a disc that intersects $D$
transversely in a single smooth point (such a loop is called a {\it geometric
generator} of $\pi_1(\CP^2-D)$) by $\theta$ is necessarily a transposition.
The smoothness of $X$ above a singular point of $D$ implies some
compatibility properties on these transpositions (geometric generators
corresponding to the two branches of $D$ at a node must map to disjoint
commuting transpositions, while to a cusp must correspond a pair of adjacent
transpositions). Finally, the connectedness of $X$ implies the surjectivity
of $\theta$ (because the subgroup $\mathrm{Im}(\theta)$ is generated by
transpositions and acts transitively on the fiber of the covering).

It must be mentioned that the amount of information present in the monodromy
morphism $\theta$ is fairly small: a classical conjecture in algebraic
geometry (Chisini's conjecture, solved by Kulikov \cite{Ku}) asserts that,
given an algebraic singular plane curve $D$ with cusps and nodes, a
symmetric group-valued monodromy morphism $\theta$ compatible with $D$
(in the above sense), if it exists, is unique except in a small finite list
of cases. Whether Chisini's conjecture also holds for symplectic branch
curves is an open question, but in any case the number of possibilities for
$\theta$ is always finite.

The study of a singular complex curve $D\subset\CP^2$ can be carried out
using the braid monodromy techniques developed in complex algebraic geometry
by Moishezon and Teicher \cite{Mo2,Te1}: the idea is to choose a linear
projection $\pi:\CP^2-\{\mathrm{pt}\}\to\CP^1$, for example
$\pi(x\!:\!y\!:\!z)=(x\!:\!y)$, in such a way that the curve $D$ lies in
general position with respect to the fibers of $\pi$, i.e.\ $D$ is
positively transverse to the fibers of $\pi$ everywhere except at isolated
non-degenerate smooth complex tangencies. The restriction $\pi_{|D}$ is then
a singular branched covering of degree $d=\deg D$, with {\it special points}
corresponding to the singularities of $D$ (nodes and cusps) and to the
tangency points. Moreover, we can assume that all special points lie in
distinct fibers of $\pi$. A plane curve satisfying these topological
requirements is called a {\it braided} (or {\it Hurwitz}) curve.
\medskip

\begin{center}
\setlength{\unitlength}{0.8mm}
\begin{picture}(80,55)(-40,-15)
\put(0,-2){\vector(0,-1){8}}
\put(2,-7){$\pi:[x\!:\!y\!:\!z]\mapsto [x\!:\!y]$}
\put(-40,-15){\line(1,0){80}}
\put(-38,-12){ $\CP^1$}
\put(-40,0){\line(1,0){80}}
\put(-40,40){\line(1,0){80}}
\put(-40,0){\line(0,1){40}}
\put(40,0){\line(0,1){40}}
\put(-38,33){ $\CP^2-\{\infty\}$}
\put(27,31){ $D$}
\multiput(-20,20)(0,-2){18}{\line(0,-1){1}}
\multiput(-5,20)(0,-2){18}{\line(0,-1){1}}
\multiput(15,15)(0,-2){9}{\line(0,-1){1}}
\multiput(15,-9)(0,-2){3}{\line(0,-1){1}}
\put(-20,-15){\circle*{1}}
\put(-5,-15){\circle*{1}}
\put(15,-15){\circle*{1}}
\qbezier[200](25,35)(5,30)(-5,20)
\qbezier[60](-5,20)(-10,15)(-15,15)
\qbezier[40](-15,15)(-20,15)(-20,20)
\qbezier[40](-20,20)(-20,25)(-15,25)
\qbezier[60](-15,25)(-10,25)(-5,20)
\qbezier[110](-5,20)(0,15)(15,15)
\qbezier[320](15,15)(5,15)(-30,5)
\put(-20,20){\circle*{1}}
\put(-5,20){\circle*{1}}
\put(15,15){\circle*{1}}
\end{picture}
\end{center}
\medskip

Except for those which contain special points of $D$, the fibers of $\pi$
are lines intersecting the curve $D$ in $d$ distinct points. If one chooses
a reference point $q_0\in\CP^1$ (and the corresponding fiber $\ell\simeq\C
\subset\CP^2$ of $\pi$), and if one restricts to an affine subset in order
to be able to trivialize the fibration $\pi$, the topology of the branched
covering $\pi_{|D}$ can be described by a {\it braid monodromy} morphism
\begin{equation}
\rho:\pi_1(\C-\{\mathrm{pts}\},q_0)\to B_d,
\end{equation}
where $B_d$ is the braid group on $d$ strings. The braid $\rho(\gamma)$
corresponds to the motion of the $d$ points of $\ell\cap D$ inside the
fibers of $\pi$ when moving along the loop $\gamma$.

Recall that the braid group $B_d$ is the fundamental group of the
configuration space of $d$ distinct points in $\R^2$; it is also the
group of isotopy classes of compactly supported orientation-preserving
diffeomorphisms of $\R^2$ leaving invariant a set of $d$ given distinct
points. It is generated by the standard {\it half-twists} $X_1,\dots,X_{d-1}$
(braids which exchange two consecutive points by rotating them 
counterclockwise by 180 degrees around each other), with relations
$X_iX_j=X_jX_i$ for $|i-j|\ge 2$ and $X_iX_{i+1}X_i=X_{i+1}X_iX_{i+1}$
(the reader is referred to Birman's book \cite{Bi} for more details).

Another equivalent way to consider the monodromy of a braided curve is
to choose an ordered system of generating
loops in the free group $\pi_1(\C-\{\mathrm{pts}\},q_0)$. The morphism
$\rho$ can then be described by a {\it factorization} in the braid group
$B_d$, i.e.\ a decomposition of the monodromy at infinity into the product
of the individual monodromies around the various special points of $D$.
By observing that the total space of $\pi$ is the line bundle $O(1)$ over
$\CP^1$, it is easy to see that the monodromy at infinity is given by the
central element $\Delta^2=(X_1\dots X_{d-1})^d$ of $B_d$ (called ``full
twist'' because it represents a rotation of a large disc by 360 degrees).
The individual monodromies around the special points are conjugated to
powers of half-twists, the exponent being $1$ in the case of tangency
points, $2$ in the case of positive nodes (or $-2$ for negative nodes), and
$3$ in the case of cusps.

The braid monodromy $\rho$ and the corresponding factorization depend on
trivialization choices, which affect them by {\it simultaneous conjugation} by
an element of $B_d$ (change of trivialization of the fiber $\ell$ of $\pi$),
or by {\it Hurwitz operations} (change of generators of the group
$\pi_1(\C-\{\mathrm{pts}\},q_0)$). There is a one-to-one correspondence between
braid monodromy morphisms $\rho:\pi_1(\C-\{\mathrm{pts}\})\to B_d$ (mapping
generators to suitable powers of half-twists) up to
these two algebraic operations and singular (not necessarily complex)
braided curves of degree $d$ in $\CP^2$ 
up to isotopy among such curves (see e.g.\ \cite{KK} for a detailed
exposition). Moreover, it is easy to check that
every braided curve in $\CP^2$ can be deformed into a braided symplectic
curve, canonically up to isotopy among symplectic braided curves (this
deformation is performed by collapsing the curve $D$ into a neighborhood
of a complex line in a way that preserves the fibers of $\pi$). However,
the curve $D$ is isotopic to a complex curve only for certain specific
choices of the morphism $\rho$.

Unlike the case of complex curves, it is not clear {\it a priori} that the
symplectic branch curve $D_k$ of one of the covering maps given by
Theorem 4.3 can be made compatible with the linear projection $\pi$;
making the curve $D_k$ braided relies on an improvement of Theorem 4.3 in
order to control more precisely the behavior of $D_k$ near its special
points (tangencies, nodes, cusps) \cite{AK}. Moreover, one must take into
account the possible occurrence of creations or cancellations of admissible
pairs of nodes in the branch curve $D_k$, which affect the braid monodromy
morphism $\rho_k:\pi_1(\C-\{\mathrm{pts}\})\to B_d$ by insertion or deletion
of pairs of factors. The uniqueness statement in Theorem 4.3 then leads to
the following result, obtained in collaboration with Katzarkov \cite{AK}:

\begin{theorem}[A.-Katzarkov]
For given large enough $k$, the monodromy morphisms $(\rho_k,\theta_k)$
associated to the approximately holomorphic branched covering maps
$f_k:X\to\CP^2$ defined by triples of sections of $L^{\otimes k}$ are, up
to conjugation, Hurwitz operations, and insertions/deletions, invariants of
the symplectic manifold $(X,\omega)$. Moreover, these invariants are {\em
complete}, in the sense that the data $(\rho_k,\theta_k)$ are sufficient to
reconstruct the manifold $(X,\omega)$ up to symplectomorphism.
\end{theorem}

It is interesting to mention that the symplectic Lefschetz pencils
constructed by Donaldson (Theorem 2.7) can be recovered very easily from
the branched covering maps $f_k$, simply by considering the $\CP^1$-valued
maps $\pi\circ f_k$. In other words, the fibers $\Sigma_{k,\alpha}$ of the
pencil are the preimages by $f_k$ of the fibers of $\pi$, and the singular
fibers of the pencil correspond to the fibers of $\pi$ through the tangency
points of $D_k$.

In fact, the monodromy morphisms $\psi_k$ of the Lefschetz pencils $\pi\circ
f_k$ (see \S 2.4) can be recovered very explicitly from $\theta_k$ and
$\rho_k$. By restriction to the line $\bar\ell=\ell\cup \{\infty\}$, the
$S_N$-valued morphism $\theta_k$ describes the topology of a fiber
$\Sigma_k$ of the pencil as an $N$-fold covering of $\CP^1$ with $d$ branch
points; the set of base points $Z_k$ is the preimage of the point at
infinity in $\bar\ell$. This makes it possible to define a {\it lifting
homomorphism} from a subgroup $B_d^0(\theta_k)\subset B_d$ ({\it liftable
braids}) to the mapping class group $\mathrm{Map}(\Sigma_k,Z_k)=
\mathrm{Map}_{g,N}$. The various monodromies are then related by the
following formula \cite{AK}:
\begin{equation}\label{eq:lifting}\psi_k=(\theta_k)_*\circ \rho_k.\end{equation}

The lifting homomorphism $(\theta_k)_*$ maps liftable half-twists to Dehn
twists, so that the tangencies between the branch curve $D_k$ and the fibers
of $\pi$ determine explicitly the vanishing cycles of the Lefschetz pencil
$\pi\circ f_k$. On the other hand, the monodromy around a node or cusp of
$D_k$ lies in the kernel of $(\theta_k)_*$.

The lifting homomorphism $\theta_*$ can be defined more precisely as
follows: the space $\tilde{\mathcal{X}}_d$ of configurations of $d$
distinct points in $\R^2$ together with branching data (a transposition
in $S_N$ attached to each point) is a finite covering of the space
$\mathcal{X}_d$ of configurations of $d$ distinct points. The morphism
$\theta$ determines a lift $\tilde{*}$ of the base point in $\mathcal{X}_d$,
and the liftable braid subgroup of $B_d=\pi_1(\mathcal{X}_d,*)$ is
the stabilizer of $\theta$ for the action of $B_d$ by deck transformations
of the covering $\tilde{\mathcal{X}}_d\to \mathcal{X}_d$, i.e.\
$B_d^0(\theta)=\pi_1(\tilde{\mathcal{X}}_d,\tilde{*})$.
Moreover, $\tilde{\mathcal{X}}_d$ is naturally equipped with a 
universal fibration $\mathcal{Y}_d\to \tilde{\mathcal{X}}_d$ by genus $g$
Riemann surfaces with $N$ marked points: the lifting homomorphism
$\theta_*:B_d^0(\theta)\to \mathrm{Map}_{g,N}$
is by definition the monodromy of this fibration.

The relation (\ref{eq:lifting}) is very useful for explicit calculations of
the monodromy of Lefschetz pencils, which is accessible to direct methods
only in a few very specific cases. By comparison, the various available
techniques for braid monodromy calculations \cite{Mo3,Te1,ADKY} are much
more powerful, and make it possible to carry out calculations in a much
larger number of cases (see \S 4.5 below). In particular, in view of
Theorem 2.7 we are mostly interested in the monodromy of high degree
Lefschetz pencils, where the fiber genus and the number of singular fibers
are very high, making them inaccessible to direct calculation even for the
simplest complex algebraic surfaces.

\subsection{Monodromy invariants in higher dimension}

When $\dim X>4$, the topology of a projective map $f_k:X-Z_k\to\CP^2$ (as
given by Theorem 4.4) and of the discriminant curve $D_k\subset\CP^2$ can
be described using the same methods as for a branched covering; the only
difference is that the morphism $\theta_k$ describing the monodromy of the
fibration above the complement of $D_k$ now takes values in the symplectic
mapping class group $\mathrm{Map}^\omega(\Sigma_k,Z_k)$ of the generic fiber
of $f_k$. Theorem 4.6 then remains true in this context \cite{Au3}.
However, we still face the same difficulty as in the case of Lefschetz
pencils of arbitrary dimension, namely the fact that the mapping
class group of a symplectic manifold of dimension $4$ or more is
essentially unknown. Therefore, even though the monodromy morphisms
$(\rho_k,\theta_k)$ provide an appealing description of a symplectic
6-manifold fibered above $\CP^2$, they cannot be used to accurately
describe a symplectic manifold of dimension $8$ or more. 
The definition of invariants using maps to higher-dimensional
projective spaces (in order to decrease the dimension of the fibers),
although possible in principle,
does not seem to be a satisfactory solution to this
problem, because the structure of the discriminant hypersurface in
$\CP^m$ becomes very complicated for $m\ge 3$.

However, one can work around this difficulty by means of a dimensional
reduction process. Indeed, the restriction of $f_k$ to the line $\bar\ell
\subset\CP^2$ defines a Lefschetz pencil structure on the symplectic
hypersurface $W_k=f_k^{-1}(\bar\ell)\subset X$, with generic fiber
$\Sigma_k$, base locus $Z_k$, and monodromy $\theta_k$.

\setlength{\unitlength}{0.68mm}
\begin{center}
\begin{picture}(120,78)(-50,-28)
\put(-40,0){\line(1,0){80}}
\put(-40,40){\line(1,0){80}}
\put(-45,-15){\line(1,0){90}}
\put(-65,-20){$\CP^2$}
\put(-55,-25){\line(1,0){90}}
\put(-55,-25){\line(1,1){20}}
\put(35,-25){\line(1,1){20}}
\put(55,-25){\line(1,1){20}}
\put(-49,34){\line(1,1){12}}
\put(-49,-6){\line(1,1){12}}
\put(38,-5.2){\line(1,1){12}}
\put(38,34.8){\line(1,1){12}}
\put(-60,20){$X$}
\put(70,-20){$\CP^1$}
\put(50,-15){\vector(1,0){10}}
\put(55,-13.5){\makebox(0,0)[cb]{$\pi$}}
\qbezier[35](-40,40)(-35,40)(-35,35)
\qbezier[35](-35,35)(-35,30)(-36,27.5)
\qbezier[35](-36,27.5)(-37,25)(-37,20)
\qbezier[35](-40,0)(-35,0)(-35,5)
\qbezier[35](-35,5)(-35,10)(-36,12.5)
\qbezier[35](-36,12.5)(-37,15)(-37,20)
\qbezier[35](-40,40)(-45,40)(-45,35)
\qbezier[35](-45,35)(-45,30)(-44,27.5)
\qbezier[35](-44,27.5)(-43,25)(-43,20)
\qbezier[35](-40,0)(-45,0)(-45,5)
\qbezier[35](-45,5)(-45,10)(-44,12.5)
\qbezier[35](-44,12.5)(-43,15)(-43,20)
\qbezier[30](-41,35)(-37,32)(-41,29)
\qbezier[20](-40,34)(-42,32)(-40,30)
\qbezier[30](-41,5)(-37,8)(-41,11)
\qbezier[20](-40,6)(-42,8)(-40,10)
\qbezier[35](40,40)(35,40)(35,35)
\qbezier[35](35,35)(35,30)(36,27.5)
\qbezier[35](36,27.5)(37,25)(37,20)
\qbezier[35](40,0)(35,0)(35,5)
\qbezier[35](35,5)(35,10)(36,12.5)
\qbezier[35](36,12.5)(37,15)(37,20)
\qbezier[35](40,40)(45,40)(45,35)
\qbezier[35](45,35)(45,30)(44,27.5)
\qbezier[35](44,27.5)(43,25)(43,20)
\qbezier[35](40,0)(45,0)(45,5)
\qbezier[35](45,5)(45,10)(44,12.5)
\qbezier[35](44,12.5)(43,15)(43,20)
\qbezier[30](39,35)(43,32)(39,29)
\qbezier[20](40,34)(38,32)(40,30)
\qbezier[30](39,5)(43,8)(39,11)
\qbezier[20](40,6)(38,8)(40,10)
\qbezier[35](-10,40)(-5,40)(-5,35)
\qbezier[35](-5,35)(-5,30)(-6,27.5)
\qbezier[35](-6,27.5)(-7,25)(-7,20)
\qbezier[35](-10,0)(-5,0)(-5,5)
\qbezier[35](-5,5)(-5,10)(-6,12.5)
\qbezier[35](-6,12.5)(-7,15)(-7,20)
\qbezier[35](-10,40)(-15,40)(-15,36)
\qbezier[35](-15,36)(-15,32)(-14,28)
\qbezier[35](-14,28)(-13,26)(-13,25)
\qbezier[35](-10,0)(-15,0)(-15,5)
\qbezier[35](-15,5)(-15,10)(-14,12.5)
\qbezier[70](-14,12.5)(-13,15)(-13,25)
\qbezier[30](-8,32)(-8,36)(-10,36)
\qbezier[30](-10,36)(-12,36)(-14,32)
\qbezier[20](-8,32)(-8,30)(-10,30)
\qbezier[20](-10,30)(-12,30)(-14,32)
\qbezier[30](-11,5)(-7,8)(-11,11)
\qbezier[20](-10,6)(-12,8)(-10,10)
\qbezier[35](10,40)(15,40)(15,35)
\qbezier[35](15,35)(15,30)(14,27.5)
\qbezier[35](14,27.5)(13,25)(13,20)
\qbezier[35](10,0)(15,0)(15,5)
\qbezier[35](15,5)(15,10)(14,12.5)
\qbezier[35](14,12.5)(13,15)(13,20)
\qbezier[35](10,0)(5,0)(5,4)
\qbezier[35](5,4)(5,8)(6,12)
\qbezier[35](6,12)(7,14)(7,15)
\qbezier[35](10,40)(5,40)(5,35)
\qbezier[35](5,35)(5,30)(6,27.5)
\qbezier[70](6,27.5)(7,25)(7,15)
\qbezier[30](12,8)(12,4)(10,4)
\qbezier[30](10,4)(8,4)(6,8)
\qbezier[20](12,8)(12,10)(10,10)
\qbezier[20](10,10)(8,10)(6,8)
\qbezier[30](9,35)(13,32)(9,29)
\qbezier[20](10,34)(8,32)(10,30)
\put(20,20){ $W$}
\put(-44,20){ $\Sigma$}
\put(-46,-23){$D$}
\put(-30,-13){$\bar\ell$}
\put(-14,32){\circle*{2}}
\put(6,8){\circle*{2}}
\qbezier[20](-44.8,32)(-42.9,31)(-41,32)
\qbezier[8](-44.8,32)(-42.9,33)(-41,32)
\multiput(-10,-1)(0,-2){7}{\line(0,-1){1}}
\multiput(10,-1)(0,-2){7}{\line(0,-1){1}}
\put(10,-15){\circle*{2}}
\put(-10,-15){\circle*{2}}
\qbezier[270](-40,-22)(-2,-15)(0,-10)
\qbezier[160](25,-18.5)(2,-15)(0,-10)
\end{picture}
\end{center}

This structure can be enriched by adding an extra section of $L^{\otimes k}$
in order to obtain a map from $W_k$ to $\CP^2$ with generic local models and
braided discriminant curve; this map can in turn be characterized using
monodromy invariants. This process can be repeated on successive hyperplane
sections until we reduce to the 4-dimensional case. Hence, given a
symplectic manifold $(X^{2n},\omega)$ and a large integer $k$, we obtain
$n-1$ singular curves $D_{k,(n)},D_{k,(n-1)},\dots,D_{k,(2)}\subset\CP^2$,
described by $n-1$ braid monodromy morphisms, and a homomorphism
$\theta_{k,(2)}$ from $\pi_1(\CP^2-D_{k,(2)})$ to a symmetric group.

These invariants are sufficient to successively reconstruct the various
submanifolds of $X$ involved in the reduction process. Indeed, given a symplectic
manifold $\Sigma_{k,(r-1)}$ of dimension $2r-2$ equipped with a Lefschetz
pencil structure with generic fiber $\Sigma_{k,(r-2)}$ and monodromy
$\theta_{k,(r)}$, and given the braid monodromy of a plane curve
$D_{k,(r)}\subset\CP^2$, we obtain a symplectic manifold $\Sigma_{k,(r)}$ of 
dimension $2r$ equipped with a $\CP^2$-valued map with generic fiber
$\Sigma_{k,(r-2)}$, monodromy $\theta_{k,(r)}$ and discriminant cuve
$D_{k,(r)}$. Composing with the projection $\pi$ from $\CP^2$ to $\CP^1$ we
obtain on $\Sigma_{k,(r)}$ a Lefschetz pencil structure with generic fiber
$\Sigma_{k,(r-1)}$, whose monodromy $\theta_{k,(r+1)}$ can be obtained from
the braid monodromy of $D_{k,(r)}$ using the lifting homomorphism induced
by $\theta_{k,(r)}$, using a formula similar to (\ref{eq:lifting}).

Therefore, the data consisting of $n-1$ braid monodromies and a symmetric
group-valued homomorphism, satisfying certain compatibility conditions
(each braid monodromy must be contained in an appropriate liftable braid
subgroup, and the monodromies around nodes and cusps must lie in the kernel
of the lifting homomorphism), completely determines a symplectic manifold
$(X,\omega)$ up to symplectomorphism. Hence, we have the following
combinatorial description of compact symplectic manifolds \cite{Au3}:

\begin{theorem}
For given $k\gg 0$, the monodromy morphisms $\rho_{k,(n)},\dots,
\rho_{k,(2)}$ with values in braid groups and $\theta_{k,(2)}$ with values
in a symmetric group that can be associated to linear systems of suitably
chosen sections of $L^{\otimes k}$ are, up to conjugations, Hurwitz
operations and insertions/deletions, invariants of the symplectic manifold
$(X^{2n},\omega)$. Moreover, these invariants contain enough information
to reconstruct the manifold $(X^{2n},\omega)$ up to symplectomorphism.
\end{theorem}

For higher-dimensional symplectic manifolds, this strategy of approach
appears to be more promising than the direct study of $\CP^m$-valued maps
for $m\ge 3$, because it makes it unnecessary to handle the very complicated
singularities present in higher-dimensional discriminant loci. However, it
must be mentioned that, whatever the approach, explicit calculations are to
this date possible only on very specific low degree examples.

\subsection{Calculation techniques}

In principle, Theorems 4.6 and 4.7 reduce the classification of compact
symplectic manifolds to purely combinatorial questions concerning braid
groups and symmetric groups, and symplectic topology seems to largely reduce
to the study of certain singular plane curves, or equivalently certain words
in braid groups.

The explicit calculation of these monodromy invariants is hard in the
general case, but is made possible for a large number of complex surfaces by
the use of ``degeneration'' techniques and of approximately holomorphic
perturbations. Hence, the invariants defined by Theorem 4.6 are explicitly
computable for various examples such as $\CP^2$, $\CP^1\times\CP^1$
\cite{Mo3}, a few complete intersections (Del Pezzo or K3 surfaces)
\cite{Ro}, the Hirzebruch surface $\mathbb{F}_1$, and all double covers
of $\CP^1\times\CP^1$ branched along connected smooth algebraic curves
(among which an infinite family of surfaces of general type) \cite{ADKY}.

The degeneration technique, developped by Moishezon and Teicher
\cite{Mo3,Te1}, starts with a projective embedding of the complex surface
$X$, and deforms the image of this embedding to a singular configuration
$X_0$ consisting of a union of planes intersecting along lines. The
discriminant curve of a projection of $X_0$ to $\CP^2$ is therefore a union
of lines; the manner in which the smoothing of $X_0$ affects this curve
can be studied explicitly, by considering a certain number of standard local
models near the various points of $X_0$ where three or more planes
intersect. This method makes it possible to handle many examples in low
degree, but in the case $k\gg 0$ that we are interested in (very positive
linear systems over a fixed manifold), the calculations can only
be carried out explicitly for very simple surfaces.

In order to proceed beyond this point, it becomes more efficient to move
outside of the algebraic framework and to consider generic approximately
holomorphic perturbations of non-generic algebraic maps; the greater
flexibility of this setup makes it possible to choose more easily computable
local models. For example, the direct calculation of the monodromy
invariants becomes possible for all linear systems of the type $\pi^*O(p,q)$
on double couvers of $\CP^1\times\CP^1$ branched along connected smooth
algebraic curves of arbitrary degree \cite{ADKY}. It also becomes possible
to obtain a general ``degree doubling'' formula, describing explicitly the
monodromy invariants associated to the linear system $L^{\otimes 2k}$ in
terms of those associated to the linear system $L^{\otimes k}$ (when $k\gg
0$), both for branched covering maps to $\CP^2$ and for 4-dimensional
Lefschetz pencils \cite{AK2}.

However, in spite of these successes, a serious obstacle restricts the
practical applications of monodromy invariants: in general, they cannot
be used efficiently to distinguish homeomorphic symplectic manifolds,
because no algorithm exists to decide whether two words in a braid group
or mapping class group are equivalent to each other via Hurwitz operations.
Even if an algorithm could be found, another difficulty is due to the large
amount of combinatorial data to be handled: on a typical interesting example,
the braid monodromy data can already consist of $\sim 10^4$ factors in
a braid group on $\sim 100$ strings for very small values of the
parameter $k$, and the amount of data grows polynomially with $k$.

Hence, even when monodromy invariants can be computed, they cannot be {\it
compared}. This theoretical limitation makes it necessary to search for
other ways to exploit monodromy data, e.g.\ by considering invariants that
contain less information than braid monodromy but are easier to use in
practice.

\subsection{Fundamental groups of branch curve complements}

Given a singular plane curve $D\subset\CP^2$, e.g.\ the branch curve of a
covering, it is natural to study the fundamental group $\pi_1(\CP^2-D)$.
The study of this group for various types of algebraic curves is a classical
subject going back to the work of Zariski, and has undergone a lot of
development in the 80's and 90's, in part thanks to the work of Moishezon
and Teicher \cite{Mo2,Mo3,Te1}. The relation to braid monodromy invariants
is a very direct one: the Zariski-van Kampen theorem provides an explicit
presentation of the group $\pi_1(\CP^2-D)$ in terms of the braid monodromy
morphism $\rho:\pi_1(\C-\{\mathrm{pts}\})\to B_d$. However, if one is
interested in the case of symplectic branch curves, it is important to
observe that the introduction or the cancellation of pairs of nodes affects
the fundamental group of the complement, so that it cannot be used directly
to define an invariant associated to a symplectic branched covering.
In the symplectic world, the fundamental group of the branch curve
complement must be replaced by a suitable quotient, the {\it stabilized
fundamental group} \cite{ADKY}.

Using the same notations as in \S 4.3, the inclusion \hbox{$i:\ell-(\ell\cap
D_k)\to \mathbb{CP}^2-D_k$} of the reference fiber of the linear projection
$\pi$ induces a surjective morphism on fundamental groups; the images of the
standard generators of the free group \hbox{$\pi_1(\ell-(\ell\cap D_k))$} and their
conjugates form a subset $\Gamma_k\subset\pi_1(\mathbb{CP}^2-D_k)$ whose
elements are called {\it geometric generators}. Recall that the images of
the geometric generators by the monodromy morphism $\theta_k$ are
transpositions in $S_N$. The creation of a pair of nodes in the curve $D_k$
amounts to quotienting $\pi_1(\CP^2-D_k)$ by a relation of the form
$[\gamma_1,\gamma_2]\sim 1$, where $\gamma_1,\gamma_2\in\Gamma_k$; however,
this creation of nodes can be carried out by deforming the branched
covering map $f_k$ only if the two transpositions $\theta_k(\gamma_1)$ and
$\theta_k(\gamma_2)$ have disjoint supports. Let $K_k$ be the normal
subgroup of $\pi_1(\CP^2-D_k)$ generated by all such commutators
$[\gamma_1,\gamma_2]$. Then we have the following result \cite{ADKY}:

\begin{theorem}[A.-D.-K.-Y.]
For given $k\gg 0$, the {\em stabilized fundamental group}
$\bar{G}_k=\pi_1(\mathbb{CP}^2-D_k)/K_k$ is an invariant of the symplectic
manifold $(X^4,\omega)$.
\end{theorem}

This invariant can be calculated explicitly for the various examples where
monodromy invariants are computable ($\CP^2$, $\CP^1\times\CP^1$, some
Del Pezzo and K3 complete intersections, Hirzebruch surface $\mathbb{F}_1$,
double covers of $\CP^1\times\CP^1$); namely, the extremely complicated
presentations given by the Zariski-van Kampen theorem in terms of braid
monodromy data can be simplified in order to obtain a manageable
description of the fundamental group of the branch curve complement.
These examples lead to various observations and conjectures.

A first remark to be made is that, for all known examples, when the
parameter $k$ is sufficiently large the stabilization operation becomes
trivial, i.e.\ geometric generators associated to disjoint transpositions
already commute in $\pi_1(\CP^2-D_k)$, so that $K_k=\{1\}$ and
$\bar{G}_k=\pi_1(\CP^2-D_k)$. For example, in the case of $X=\CP^2$ with
its standard K\"ahler form, we have $\bar{G}_k=\pi_1(\CP^2-D_k)$ for all
$k\ge 3$. Therefore, when $k\gg 0$ no information seems to be lost when
quotienting by $K_k$ (the situation for small values of $k$ is very
different).

The following general structure result can be proved for the groups
$\bar{G}_k$ (and hence for $\pi_1(\CP^2-D_k)$) \cite{ADKY}:

\begin{theorem}[A.-D.-K.-Y.]
Let $f:(X,\omega)\to \CP^2$ be a symplectic branched covering of degree $N$,
with braided branch curve $D$ of degree $d$, and let
$\bar{G}=\pi_1(\CP^2-D)/K$ be the stabilized fundamental group of the branch
curve complement. Then there exists a natural exact sequence
$$1\longrightarrow G^0\longrightarrow \bar{G} \longrightarrow
S_N\times \Z_d\longrightarrow \Z_2\longrightarrow 1.$$
Moreover, if $X$ is simply connected
then there exists a natural surjective homomorphism
$\phi:G^0\twoheadrightarrow (\Z^2/\Lambda)^{n-1}$, where
$$\Lambda=\{(c_1(K_X)\cdot \alpha,[f^{-1}(\bar\ell)]\cdot \alpha),
\ \alpha\in H_2(X,\Z)\}.$$
\end{theorem}

In this statement, the two components of the morphism $\bar{G}\to S_N\times
\Z_d$ are respectively the monodromy of the branched covering,
$\theta:\pi_1(\CP^2-D)\to S_N$, and the linking number (or abelianization,
when $D$ is irreducible)
morphism $$\delta:\pi_1(\CP^2-D)\to \Z_d\ (\simeq H_1(\CP^2-D,\Z)).$$
The subgroup $\Lambda$ of $\Z^2$ is entirely determined by the numerical
properties of the canonical class $c_1(K_X)$ and of the hyperplane class
(the homology class of the preimage of a line $\bar\ell\subset\CP^2$: in
the case of the covering maps of Theorem 4.3 we have 
$[f^{-1}(\bar\ell)]=c_1(L^{\otimes k})=\frac{k}{2\pi}[\omega]$). The morphism
$\phi$ is defined by considering the $N$ lifts in $X$ of a closed loop
$\gamma$ belonging to $G^0$, or more precisely their homology classes (whose
sum is trivial) in the complement of a hyperplane section and of the
ramification curve in $X$.

Moreover, in the known examples we have a much stronger result on
the structure of the subgroups $G_k^0$ for the branch
curves of large degree covering maps (determined by sufficiently ample
linear systems) \cite{ADKY}.

Say that the simply connected complex surface $(X,\omega)$ belongs to the
class $(\mathcal{C})$ if it belongs to the list of computable examples:
$\CP^1\times\CP^1$, $\CP^2$, the Hirzebruch surface $\mathbb{F}_1$ (equipped
with any K\"ahler form), a Del Pezzo or K3 surface (equipped with a K\"ahler
form coming from a specific complete intersection realization), or a double
cover of $\CP^1\times \CP^1$ branched along a connected smooth algebraic
curve (equipped with a K\"ahler form in the class $\pi^*O(p,q)$ for $p,q\ge
1$). Then we have:

\begin{theorem}[A.-D.-K.-Y.]
If $(X,\omega)$ belongs to the class $(\mathcal{C})$, then
for all large enough $k$ the homomorphism
$\phi_k$ induces an isomorphism on the abelianized groups, i.e.\ $\mathrm{Ab}
\,G_k^0\simeq$ $(\Z^2/\Lambda_k)^{N_k-1}$, while $\mathrm{Ker}\,\phi_k=[G_k^0,
G_k^0]$ is a quotient of $\Z_2\times\Z_2$.
\end{theorem}

It is natural to make the following conjecture:

\begin{conj}
If $X$ is a simply connected symplectic 4-manifold, then for all large
enough $k$ the homomorphism
$\phi_k$ induces an isomorphism on the abelianized groups, i.e.\
$\mathrm{Ab} \,G_k^0\simeq$ $(\Z^2/\Lambda_k)^{N_k-1}$.
\end{conj}

\subsection{Symplectic isotopy and non-isotopy}

While it has been well-known for many years that compact symplectic
4-manifolds do not always admit K\"ahler structures, it has been discovered
more recently that symplectic curves (smooth or singular) in a given
manifold can also offer a wider range of possibilities than complex curves.
Proposition 4.2 and Theorem 4.3 establish a bridge between these two
phenomena: indeed, a covering of $\CP^2$ (or any other complex surface)
branched along a complex curve automatically inherits a complex structure.
Therefore, starting with a non-K\"ahler symplectic manifold, Theorem 4.3
always yields branch curves that are not isotopic to any complex curve in
$\CP^2$. The study of isotopy and non-isotopy phenomena for curves is
therefore of major interest for our understanding of the topology of
symplectic 4-manifolds.

The {\it symplectic isotopy problem} asks whether, in a given complex
surface, every symplectic submanifold representing a given homology class
is isotopic to a complex submanifold. The first positive result in this
direction was due to Gromov, who showed using his compactness result for
pseudo-holomorphic curves (Theorem 1.10) that, in $\CP^2$, a smooth
symplectic curve of degree $1$ or $2$ is always isotopic to a complex curve.
Successive improvements of this technique have made it possible to extend
this result to curves of higher degree in $\CP^2$ or $\CP^1\times\CP^1$; the
currently best known result is due to Siebert and Tian, and makes it
possible to handle the case of smooth curves in $\CP^2$ up to degree 17
\cite{ST}. Isotopy results are also known for sufficiently simple singular
curves (Barraud, Shevchishin \cite{Sh}, \dots).

Contrarily to the above examples, the general answer to the symplectic
isotopy problem appears to be negative. The first counterexamples
among smooth connected symplectic curves were found by Fintushel and Stern
\cite{FS2}, who constructed by a {\it braiding} process infinite families of
mutually non-isotopic symplectic curves representing a same homology class
(a multiple of the fiber) in elliptic surfaces, and to Smith, who used the
same construction in higher genus. However, these two constructions are
preceded by a result of Moishezon \cite{Mo4}, who established in the
early 90's a result implying the existence in $\CP^2$ of infinite families
of pairwise non-isotopic singular symplectic curves of given degree with
given numbers of node and cusp singularities. A reformulation of Moishezon's
construction makes it possible to see that it also relies on braiding;
moreover, the braiding construction can be related to a surgery operation
along a Lagrangian torus in a symplectic 4-manifold, known as {\it Luttinger
surgery} \cite{ADK}. This reformulation makes it possible to vastly simplify
Moishezon's argument, which was based on lengthy and delicate calculations
of fundamental groups of curve complements, while relating it with various
constructions developed in 4-dimensional topology.

Given an embedded Lagrangian torus $T$ in a symplectic 4-manifold
$(X,\omega)$, a homotopically non-trivial embedded loop $\gamma\subset T$
and an integer $k$, Luttinger surgery is an operation that consists in
cutting out from $X$ a tubular neighborhood of $T$, foliated by parallel
Lagrangian tori, and gluing it back in such a way that the new meridian loop
differs from the old one by $k$ twists along the loop $\gamma$ (while
longitudes are not affected), yielding a new symplectic manifold
$(\tilde{X},\tilde{\omega})$. This relatively little-known construction,
which e.g.\ makes it possible to turn a product $T^2\times \Sigma$ into any
surface bundle over $T^2$, or to transform an untwisted fiber sum into a
twisted one, can be used to described in a unified manner numerous examples
of exotic symplectic 4-manifolds constructed in the past few years.

Meanwhile, the braiding construction of symplectic curves starts with a
(possibly singular) symplectic curve $\Sigma\subset(Y^4,\omega_Y)$ and two
symplectic cylinders embedded in $\Sigma$, joined by a
Lagrangian annulus contained in the complement of $\Sigma$,
and consists in performing $k$ half-twists between these two cylinders in
order to obtain a new symplectic curve $\tilde\Sigma$ in $Y$.

\begin{center}
\epsfig{file=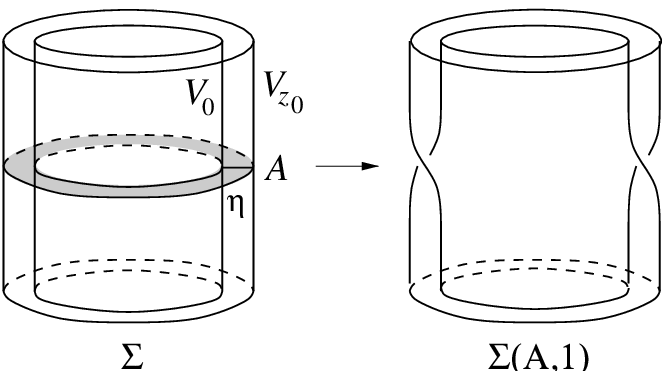,height=1.2in}
\end{center}

When $\Sigma$ is the branch curve of a symplectic branched covering $f:X\to
Y$, the following result holds \cite{ADK}:

\begin{proposition}
The covering of $Y$ branched along the symplectic curve
$\tilde\Sigma$ obtained by braiding $\Sigma$ along a Lagrangian annulus
$A\subset Y-\Sigma$ is naturally symplectomorphic to the manifold
$\tilde{X}$ obtained from the branched cover $X$ by Luttinger surgery along
a Lagrangian torus $T\subset X$ formed by the union of two lifts of $A$.
\end{proposition}

Hence, once an infinite family of symplectic curves has been constructed by
braiding, it is sufficient to find invariants that distinguish the
corresponding branched covers in order to conclude that the curves are not
isotopic. In the Fintushel-Stern examples, the branched covers are
distinguished by their Seiberg-Witten invariants, whose behavior is well
understood in the case of elliptic fibrations and their surgeries.

In the case of Moishezon's examples, a braiding construction makes it
possible to construct, starting from complex curves
$\Sigma_{p,0}\subset\CP^2$ ($p\ge 2$) of degree $d_p=9p(p-1)$ with
$\kappa_p=27(p-1)(4p-5)$ cusps and $\nu_p=27(p-1)(p-2)(3p^2+3p-8)/2$ nodes,
symplectic curves $\Sigma_{p,k}\subset \CP^2$ for all $k\in\Z$, with the same
degree and numbers of singular points. By Proposition 4.12, these curves can
be viewed as the branch curves of symplectic coverings whose total spaces
$X_{p,k}$ differ by Luttinger surgeries along a Lagrangian torus $T\subset
X_{p,0}$. The effect of these surgeries on the canonical class and on the
symplectic form can be described explicitly, which makes it possible to
distinguish the manifolds $X_{p,k}$: the canonical class of
$(X_{p,k},\omega_{p,k})$ is given by
$p\,c_1(K_{p,k})=(6p-9)[\omega_{p,k}]+(2p-3)k\,PD([T])$. Moreover, 
$[T]\in H_2(X_{p,k},\Z)$ is not a torsion class, and if $p\not\equiv 0\mod
3$ or $k\equiv 0\mod 3$ then it is a primitive class \cite{ADK}. This implies
that infinitely many of the curves $\Sigma_{p,k}$ are pairwise non-isotopic.

It is to be observed that the argument used by Moishezon to distinguish the
curves $\Sigma_{p,k}$, which relies on a calculation of the
fundamental groups $\pi_1(\CP^2-\Sigma_{p,k})$ \cite{Mo4}, is related to
the one in \cite{ADK} by means of Conjecture 4.11, of which it can be
concluded a posteriori that it is satisfied by the given branched covers
$X_{p,k}\to \CP^2$: in particular, the fact that $\pi_1(\CP^2-\Sigma_{p,k})$
is infinite for $k=0$ and finite for $k\neq 0$ is consistent with the
observation that the canonical class of $X_{p,k}$ is proportional to its
symplectic class iff $k=0$.

%}

\end{document}